\documentclass{amsart}
\usepackage[alphabetic]{amsrefs}
\usepackage[all]{xy}

\numberwithin{equation}{subsection}

\swapnumbers
\theoremstyle{plain}
\newtheorem{thm}[subsection]{Theorem}
\newtheorem{prop}[subsection]{Proposition}

\newtheorem*{thm*}{Theorem}
\newtheorem{lemma}[subsubsection]{Lemma}
\theoremstyle{definition}

\theoremstyle{remark}
\newtheorem{rem}[subsection]{Remark}
\newtheorem{rems}[subsection]{Remarks}

\renewcommand{\O}{\mathcal{O}}
\newcommand{\FF}{\mathcal{F}}

\newcommand{\Fp}{{\mathbb{F}_p}}

\newcommand{\Fq}{{\mathbb{F}_q}}

\newcommand{\Z}{\mathbb{Z}}
\newcommand{\Q}{\mathbb{Q}}
\newcommand{\R}{\mathbb{R}}
\newcommand{\C}{\mathbb{C}}
\newcommand{\A}{\mathbb{A}}
\renewcommand{\P}{\mathbb{P}}

\newcommand{\AF}{\mathcal{A}}
\newcommand{\sha}{{\hbox to 10pt{\rlap{\hskip2.8pt\vrule
height6pt\hskip1.6pt\vrule height6pt\hskip1.6pt
\vrule height6pt}\hskip1pt\vrule height0.8pt width 8pt\hskip1pt}}}

\newcommand{\tensor}{\otimes}

\newcommand{\nodiv}{\not|}
\def\nodiv{\mathrel{\mathchoice{\not|}{\not|}{\kern-.2em\not\kern.2em|}{\kern-.2em\not\kern.2em|}}}

\DeclareMathOperator{\tr}{Tr}

\DeclareMathOperator{\ord}{ord}
\DeclareMathOperator{\rk}{Rank}
\DeclareMathOperator{\dvsr}{div}

\DeclareMathOperator{\Hom}{Hom}
\DeclareMathOperator{\aut}{Aut}
\DeclareMathOperator{\pic}{Pic}
\DeclareMathOperator{\picvar}{PicVar}

\DeclareMathOperator{\gal}{Gal}
\DeclareMathOperator{\spec}{Spec}
\DeclareMathOperator{\en}{End}

\DeclareMathOperator{\divr}{\rm Div}

\newcommand{\Fpbar}{{\overline{\mathbb{F}}_p}}
\newcommand{\Fqbar}{{\overline{\mathbb{F}}_q}}

\newcommand{\kbar}{{\overline{k}}}
\newcommand{\ratto}{{\dashrightarrow}}
\newcommand{\Ql}{{\mathbb{Q}_\ell}}

\DeclareMathOperator{\Fr}{Fr}

\newcommand{\CC}{\mathcal{C}}
\newcommand{\DD}{\mathcal{D}}
\renewcommand{\SS}{\mathcal{S}}
\newcommand{\XX}{\mathcal{X}}
\newcommand{\YY}{\mathcal{Y}}

\def\clap#1{\hbox to 0pt{\hss#1\hss}} 
\def\mathllap{\mathpalette\mathllapinternal}

\def\mathllapinternal#1#2{% 
\llap{$\mathsurround=0pt#1{#2}$}}

\theoremstyle{plain}

\begin{document}
\title[Mordell-Weil groups over function fields]{
%Homomorphisms of Jacobians
%  over $k$ and Mordell-Weil groups over $k(t)$
On Mordell-Weil groups of Jacobians over function fields} 
\author{Douglas Ulmer}
\address{Department of Mathematics \\ University of Arizona \\ Tucson,
  AZ 85721}
\curraddr{School of Mathematics \\ Georgia Institute of Technology
  \\ Atlanta, GA 30332}
\email{ulmer@math.gatech.edu}
%\thanks{This paper is based upon work supported by the National
%Science Foundation under Grant No. DMS 0400877}
\date{February 2, 2011}
\begin{abstract}
  We study the arithmetic of abelian varieties over $K=k(t)$ where $k$
  is an arbitrary field.  The main result relates Mordell-Weil groups
  of certain Jacobians over $K$ to homomorphisms of other Jacobians
  over $k$.  Our methods also yield completely explicit points on
  elliptic curves with unbounded rank over $\Fpbar(t)$ and a new
  construction of elliptic curves with moderately high rank over
  $\C(t)$.
\end{abstract}
\maketitle

\section{Introduction}
In \cite{UlmerR2} we showed that large analytic ranks are ubiquitous
in towers of function fields.  For example, if $k$ is a finite field
and $E$ is an elliptic curve over $K=k(t)$ whose $j$-invariant is not
in $k$, then after replacing $k$ by a finite extension, there are
extensions of $K$ isomorphic to $k(u)$ such that the order of
vanishing at $s=1$ of the $L$-functions $L(E/k(u),s)$ is unbounded.
In many cases, the finite extension of $k$ is not needed and we may
take extensions of the form $k(t^{1/d})$ for $d$ varying through
integers prime to the characteristic of $k$.

We were also able to show that the Birch and Swinnerton-Dyer
conjecture holds for certain elliptic curves and higher genus
Jacobians and thus exhibit abelian varieties with Mordell-Weil groups
of arbitrarily large rank.  However, the class of abelian varieties
for which we can prove the BSD conjecture is much more limited than
the class for which results on large analytic rank apply.

In her thesis (see \cite{Berger08}), Lisa Berger gave a construction
of a large class of curves for which the BSD conjecture holds {\it a
  priori\/}.  (Here by {\it a priori\/} BSD, we mean that the equality
$\rk=\ord$ is shown to hold without necessarily calculating $\rk$ or
$\ord$.)  Combined with the $L$-function results of \cite{UlmerR2},
this yields {\it families\/} of abelian varieties over function fields
with large analytic and algebraic rank.  The key point in Berger's
work is a construction of towers of surfaces each of which is
dominated by a product of curves in the sense of \cite{Schoen96}.
Roughly speaking, Berger shows that certain covers of products of
curves are themselves dominated by products of curves by using a
fundamental group argument.

The aim of this article is to elucidate the geometry of Berger's
construction and to use it to study the arithmetic of abelian
varieties over $K=k(t)$ where $k$ is an arbitrary field.  The main
result relates Mordell-Weil groups of certain Jacobians over $K$ to
homomorphisms of other Jacobians over $k$.

More precisely, let $\CC$ and $\DD$ be two smooth, irreducible curves
over $k$ and choose non-constant, separable rational functions $f$ on
$\CC$ and $g$ on $\DD$.  Consider the closed subset $Y$ of
$\CC\times_k\DD\times_k\spec K$ defined by $f-tg=0$.  Under mild
hypotheses on $f$ and $g$, $Y$ is an irreducible curve and has a
smooth proper model $X\to\spec K$.  Let $J$ be the Jacobian of $X$.
The main result (Theorem~\ref{thm:main}) is a formula relating the
rank of the Mordell-Weil group of $J$ over the fields
$K_d=\kbar(t^{1/d})$ to homomorphisms between the Jacobians of covers
of $\CC$ and $\DD$ over $\kbar$.  (Here $\kbar$ is the algebraic
closure of $k$.)  Thus hard questions about certain abelian varieties
over $K_d$ are reduced to questions over the simpler field $\kbar$.

The idea of using domination by a product of curves to study
N\'eron-Severi and Mordell-Weil groups appeared already in
\cite{Shioda86}, and he used these ideas in \cite{Shioda91} to produce
explicit points on an elliptic curve of large rank over $\Fqbar(t)$.
The dominating curves in Shioda's work are quotients of Fermat curves
and the elliptic curves he studies are isotrivial.  The principal
innovation in this paper is that we systematically use families of
curves as in \cite{Berger08} and this leads to results about
non-isotrivial abelian varieties and also to interesting
specialization phenomena.  It also demonstrates that high ranks are
not fundamentally a byproduct of supersingularity.

Here is an outline of the paper: In sections~\ref{s:CandS} and
\ref{s:DPC} we discuss some background material.  In
sections~\ref{s:Berger}, \ref{s:eDPC}, and \ref{s:Formula} we explain
the construction in more detail and prove the rank formula.  The last
three sections of the paper give two additional applications of these
ideas.  One yields completely explicit points on elliptic curves with
unbounded rank over $\Fpbar(t)$ and the other is a new construction of
elliptic curves over $\C(t)$ with moderately high rank.

It is a pleasure to acknowledge the support of several institutions in
France during a sabbatical year when much of the work on this paper
was carried out: IHES, the Universit\'e de Paris-Sud, the Universit\'e
de Paris-VI, and the Universit\'e de Paris-VII.  Thanks are also due
to Kirti Joshi and Dinesh Thakur for many valuable conversations
related to this project, to Lisa Berger and Tommy Occhipinti for their
comments and questions, and to Johan de Jong and Rutger Noot for
comments on their work on CM Jacobians.  This research was partially
supported by NSF grant DMS 0701053.

\section{Curves and surfaces}\label{s:CandS}

\subsection{Definitions}
If $F$ is a field and $X$ is a scheme over $\spec F$, the phrase ``$X$ is
smooth over $F$'' means that the morphism $X\to\spec F$ is smooth.
This implies that $X$ is a regular scheme.  When $F$ is not perfect,
the converse is false: $X$ being a regular scheme does not
imply that $X$ is smooth over $F$.

If $F$ is perfect and $X$ is a noetherian scheme over $F$, then $X$ is
smooth over $F$ if and only if $X$ is regular and in this case,
$X\times_F\spec F'\to\spec F'$ is smooth for every field extension
$F'$ of $F$.  (See for example \cite{Liu}*{4.3} for this and the
previous paragraph.  Also, a particularly clear treatment appears in
Section~V.4 of a forthcoming book of Mumford and Oda.)

Throughout the paper, $k$ is a field, $\kbar$ is an algebraic closure
of $k$, and $K$ is the rational function field $k(t)$.  Note that if
$k$ has positive characteristic, then $K$ is not perfect.

A {\it variety\/} is by definition a separated, reduced scheme of
finite type over a field.  A {\it curve\/} (resp.~{\it surface\/}) is
by definition a purely 1-dimensional (resp.~2-dimensional)
variety.

\subsection{Curves and surfaces}\label{ss:CandS}
Let $\XX$ be a smooth and proper surface over $k$, equipped with a
non-constant $k$-morphism $\pi:\XX\to\P^1_k$ to the projective line
over $k$.  Let $X\to\spec K$ be the generic fiber of $\pi$, so $X$ is
a curve over $K$.  We assume that $\pi$ is generically smooth in the
sense that $X\to\spec K$ is smooth.  (This is equivalent to the
existence of a non-empty open subset $U\subset\P^1_k$ over which $\pi$
is smooth.)  If $\XX$ is proper over $k$, then $X$ is proper over $K$.

Conversely, if $X$ is a smooth curve over $K$, then there exists a
regular surface $\XX$ over $k$ equipped with a $k$-morphism
$\pi:\XX\to\P^1_k$ whose generic fiber is the given $X$.  If $X$ is
proper over $K$ then we may choose $\XX$ to be proper over $k$ and if
$X$ has genus $>0$ and we further insist that $\pi$ be relatively
minimal (i.e., have no $(-1)$-curves in its fibers), then the pair
$(\XX,\pi)$ is unique up to isomorphism.   (A convenient, modern
reference for this and related facts is \cite{Liu}*{Chapter~9}.)

In the situation of the previous paragraph, when $k$ is not perfect it
could happen that $\XX\to k$ is not smooth.  To rule this out, we
explicitly assume that $\XX\to k$ is smooth.  In the situations
considered in this paper, the starting point will always be a surface
$\XX$ smooth over $k$, so this will not be an issue.

\subsection{Points and sections}
With notation as above, there is a bijection between sections of $\pi$
and $K$-rational points of $X$.  More generally, the map
$$D\mapsto D\times_{\P^1_k}\spec K$$ 
induces a homomorphism $\divr(\XX)\to\divr(X)$ from divisors on $\XX$
to divisors on $X$.  Under this homomorphism, multi-sections (reduced
and irreducible subschemes of codimension 1 in $\XX$ flat over
$\P^1_k$) map to closed points of $X$; their residue fields are finite
extensions of $K$ of degree equal to the degree of the multisection
over $\P^1_k$.  Divisors on $\XX$ supported in fibers of $\pi$ maps to
the zero (empty) divisor on $X$.  The homomorphism
$\divr(\XX)\to\divr(X)$ is surjective since a closed point of $X$ is
the generic fiber of its scheme-theoretic closure in $\XX$.

\subsection{Shioda-Tate}\label{ss:S-T}
In this subsection, we review the Shioda-Tate formula
(\cite{Shioda99}), which relates the N\'eron-Severi group of $\XX$ to
the Mordell-weil group of the Jacobian of $X$.  Throughout, we assume
that $k$ is algebraically closed.

Assume that $\XX$ is a geometrically irreducible, smooth, proper
surface over $k$ with a $k$-morphism $\pi:\XX\to\P^1_k$ whose generic
fiber $X\to \spec K$ is smooth.

Let $\divr(\XX)$ be the group of divisors on $\XX$, let $\pic(\XX)$ be
the Picard group of $\XX$, and let $NS(\XX)$ be the N\'eron-Severi
group of $\XX$.  
%By definition, $NS(\XX)$ is the image of $\pic(\XX)$
%in $NS(\XX\times_k\overline{k})$.  
Let $\pic^0(\XX)$ be the kernel of the natural surjection
$\pic(\XX)\to NS(\XX)$.

%[Kleiman 9.5.10 shows that this gives the connected component of $\pic$
%as $\pic^0$]

Let $\divr(X)$ be the group of divisors on $X$, let $\pic(X)$ be the
Picard group of $X$, $\pic^0(X)$ its degree zero subgroup, and $J=J_X$
the Jacobian variety of $X$, so that $\pic^0(X)=J(K)$.

Recall the natural homomorphism $\divr(\XX)\to\divr(X)$ induced by
taking the fiber product with $K$.  For a divisor $D\in\divr(\XX)$,
write $D.X$ for the degree of the image of $D$ in $\divr(X)$.  Let
$L^1\divr(\XX)$ be the subgroup of divisors $D$ such that $D.X=0$ and
let $L^1\pic(\XX)\subset\pic(\XX)$ and $L^1NS(\XX)\subset NS(\XX)$ be
the subgroups generated by $L^1\divr(\XX)$.  It is clear that
$NS(\XX)/L^1NS(\XX)$ is an infinite cyclic group.  The homomorphism
$L^1\divr(\XX)\to\divr^0(X)$ (divisors of degree 0) descends to the
Picard group and yields an exact sequence
\begin{equation}\label{eq:restriction}
0\to L^2\pic(\XX)\to L^1\pic(\XX)\to\pic^0(X)\to0
\end{equation}
where $L^2\pic(\XX)$ is defined to make the sequence exact.  A simple
geometric argument (see p.~363 of \cite{Shioda99}) shows that
$L^2\pic(\XX)$ is the subgroup of $\pic(\XX)$ generated by the classes
of divisors supported in the fibers of $\pi$.  

Let $L^2NS(\XX)$ be the image of $L^2\pic(\XX)$ in $NS(\XX)$.  It is
well known that $L^2NS(\XX)$ is free abelian of rank $1+\sum_v(f_v-1)$
where the sum is over closed points of $\P^1$ and $f_v$ is the number
of irreducible components of the fiber over $v$.  (The main point
being the fact that the intersection form restricted to the components
of a fiber of $\pi$ is negative semi-definite, with kernel equal to
rational multiples of the fiber.  See, e.g.,
\cite{Beauville}*{VIII.4}.)

Let $\picvar(\XX)$ be the Picard variety of $\XX$, so that
$\pic^0(\XX)=\picvar(\XX)(k)$.  The homomorphism
$L^1\divr(\XX)\to\divr^0(X)$ induces a morphism of abelian varieties
$\picvar(\XX)\to J$.  Let $(B,\tau)$ be the $K/k$-trace of $J$.  (See
\cite{Conrad06} for a modern account.)  By the universal property of
the $K/k$ trace, the homomorphism above factors through
$\picvar(\XX)\to B$.

\begin{prop}\label{prop:ST}
  \textup{(Shioda-Tate)}  With notations and hypotheses as above, we have:
\begin{enumerate}
\item The canonical morphism of $k$-abelian varieties
  $\picvar(\XX)\to B=\tr_{K/k}J$ is an isomorphism.
\item Defining the Mordell-Weil group of $J$ to be $MW(J)=J(K)/\tau
  B(k)$, we have
$$MW(J)\cong\frac{L^1NS(\XX)}{L^2NS(\XX)}.$$
\item In particular,
$$\rk MW(J)=\rk NS(\XX)-2-\sum_v(f_v-1)$$
where the sum is over the closed points of $\P^1_k$ and $f_v$ is the
number of irreducible components in the fiber of $\pi$ over $v$.
\end{enumerate}
\end{prop}

It is well known, but apparently not written down in detail, that this
result holds under considerably more general hypotheses, for example
when $k$ is finite.  We will discuss this and related questions in a
future publication.

\subsection{Heights}\label{ss:heights}
In \cite{Shioda99}, Shioda also defines a height pairing on $MW(J)$.
We briefly review the construction.

Choose a class $c\in MW(J)$.  Shioda proves that there is an element
$[\tilde D]\in L^1NS(\XX)\tensor\Q$ which is congruent to $c$ in
$MW(J)\tensor\Q \cong(L^1NS(\XX)/L^2NS(\XX))\tensor\Q$ and which is
orthogonal to $L^2(\XX)$ under the intersection pairing.  Roughly
speaking, one ``lifts'' $c$ to a divisor $D$ on $\XX$ whose class is
in $L^1(\XX)$ and then ``corrects'' $D$ by an element of
$L^2NS(\XX)\tensor\Q$ to ensure orthogonality with $L^2NS(\XX)$.  The
resulting $[\tilde D]$ is unique up to the addition of a multiple of
$F$, the class of a fiber of $\pi$ in $NS(\XX)$.  If $c$ and $c'$ are
elements of $MW(J)$ with corresponding divisors $\tilde D$ and $\tilde
D'$, then Shioda defines their inner product by
$$\langle c,c'\rangle = -{\tilde D}.{\tilde D'}$$
where the dot on the right signifies the intersection product on
$\XX$.  Shioda proves that $\langle\cdot,\cdot\rangle$ yields a
positive definite symmetric inner product on the real vector space
$MW(J)\tensor \R$.

When $k=\Fq$ is finite, it is well known (and not difficult to show)
that $\langle\cdot,\cdot\rangle\log q$ enjoys the properties
characterizing the N\'eron-Tate canonical height.

\section{Domination by a product of curves}\label{s:DPC}

\subsection{DPC}
A variety $\XX$ over $k$ is said to be ``dominated by a product of
curves'' (or ``DPC'') if there is a dominant rational map
$\prod_i\CC_i\ratto\XX$ where each $\CC_i$ is a curve over $k$.  This
notion was introduced and studied by Schoen in \cite{Schoen96},
following Deligne \cite{Deligne72}.  If $\XX$ has dimension $n$ and is
DPC then it admits a dominant rational map from a product of $n$
curves (\cite{Schoen96}*{6.1}).  We will be concerned only with the
case $n=2$, and therefore with rational maps from products of two
curves $\CC\times_k \DD$ to a surface $\XX$.

Schoen develops several invariants related to the notion of DPC and
obstructions to a variety being DPC.  In particular, he shows that for
any field $k$, there are surfaces over $k$ which are not DPC.

\subsection{DPC and NS}
Part of the interest of the notion of DPC is that the N\'eron-Severi
group of a product of curves is well-understood.  Indeed, (under the
assumption that $\CC$ and $\DD$ have $k$-rational points), we have a
canonical isomorphism
$$NS(\CC\times\DD)\cong\Z^2\times\Hom_{k-av}(J_\CC,J_\DD)$$
where $J_\CC$ and $J_\DD$ are the Jacobians of $\CC$ and $\DD$ and the
homomorphisms are of abelian varieties over $k$.  The factor $\Z^2$
corresponds to classes of divisors which are ``vertical'' or
``horizontal'', i.e., contained in the fibers of the projections to
$\CC$ or $\DD$.  The map from $NS$ to $\Hom$ sends a divisor
class to the action of the induced correspondence on the Jacobians.
(This is classical.  See \cite{MilneJacobians}*{Cor.~6.3} for a
modern account.)

If $\CC\times\DD\ratto\XX$ is a dominant rational map, then after
finitely many blow-ups, we have a morphism
$\widetilde{\CC\times\DD}\to\XX$.  The effect of a blow-up on $NS$ is
well-known (each blow-up adds a factor of $\Z$ to $NS$) and so we have
good control of $NS(\widetilde{\CC\times\DD})$.  In favorable cases,
then can be parlayed into good control of $NS(\XX)$.

\subsection{Finitely generated $k$}
In the case where $k$ is finitely generated over its prime field,
there is a further favorable consequence of DPC, namely the Tate
conjecture on divisors and cohomology classes.  Indeed, this
conjecture has been proven for products of curves (by Tate
when $k$ is finite \cite{Tate66}, by Zarhin when $k$ is finitely
generated of characteristic $p>0$ \cite{Zarhin74}, and by Faltings
when $k$ is finitely generated of characteristic 0 \cite{Faltings83}).
It follows immediately for products of curves blown up at finitely
many closed points and it then descends to images of these varieties.
The details of the descent are explained in \cite{TateMotives}.  The
upshot is that for $\XX$ a smooth and proper variety of any dimension
over a finitely generated field $k$ which is DPC, and for any $\ell$
not equal to the characteristic of $k$, we have
$$\rk NS(\XX)=
\dim_{\Ql}H^2(\XX\times\overline{k},\Ql)(1)^{\gal(\overline{k}/k)}.$$ The
cohomology on the right is $\ell$-adic \'etale cohomology and the
superscript denotes the space of invariants under the natural action
of Galois.

\subsection{Finite $k$}
Now assume that $k$ is finite, $\XX$ is a smooth and proper surface
over $k$, and $\XX$ is equipped with a morphism $\pi:\XX\to\P^1$ as in
Subsection~\ref{ss:CandS}.  Let $X/K$ be the generic fiber of $\pi$
and $J=J_X$ the Jacobian of $X$.  If $\XX$ is DPC, then the Tate
conjecture holds for $\XX$.  It is well known that the Tate conjecture
for $\XX$ is equivalent to the Birch and Swinnerton-Dyer (BSD)
conjecture for $J$ (\cite{Tate68}).  Thus DPC of $\XX$ gives {\it a
  priori} BSD for $J$.  Explicitly, we have
$$\rk MW(J)=\ord_{s=1}L(J,s)=\ord_{s=1}L(X,s).$$
(We also have the finiteness of $\sha(J/K)$ and the refined conjecture
of Birch and Swinnerton-Dyer, but we will not discuss these facts in
detail here.)

\subsection{DPCT}
As noted in the introduction, when $k$ is finite \cite{UlmerR2} gives
a large supply of curves $X$ over $K$ such that $\ord_{s=1}L(X/K_d,s)$
is unbounded as $d$ varies through integers prime to the
characteristic of $k$ and $K_d=k(t^{1/d})$.  If one has the BSD
conjecture for $J=Jac(X)$ in each layer of the tower, then one obtains
many examples of Mordell-Weil groups of unbounded rank.

As we have seen, domination by a product of curves provides {\it a
  priori\/} BSD.  However, we arrive at the following problem: suppose
that $\XX$ is DPC and $\pi:\XX\to\P^1$ is as in
Subsection~\ref{ss:CandS}, so that we have BSD for $J=J_X$ over $K$.
Let $\XX_d$ be a smooth proper model over $k$ of the fiber product
\begin{equation*}
\xymatrix{\P^1_k\times_{\P^1_k}\XX\ar[r]\ar[d]&\XX\ar[d]^{\pi_1}\\
\P^1_k\ar[r]_{r_d}&\P^1_k}
\end{equation*}
where $r_d^*(t)=u^d$.  Then it is not true in general that $\XX_d$ is
DPC and we cannot conclude that BSD holds for $J_X/K_d$.  The
motivation for Berger's thesis was to find a large supply of surfaces
$\XX$ which are fibered over $\P^1$ such that the base changed
surfaces as above remain DPC.  She calls this ``domination by a
product of curves in a tower'' or ``DPCT''.  If $\XX\to\P^1$ is DPCT,
then the generic fibers of $\pi_d:\XX_d\to\P^1$ are curves over
$K_d=k(u)=k(t^{1/d})$ which will often have unbounded analytic and
algebraic rank.  In fact, Berger's construction leads to families of
elliptic curves over $K$ with unbounded rank.

The aim of the rest of this paper is to elucidate the geometry of
Berger's construction and to use it to study ranks over arbitrary
fields $k$.

\section{Berger's construction}\label{s:Berger}

\subsection{Special pencils}\label{ss:SpecialPencils}
Fix two smooth, proper, irreducible curves $\CC$ and $\DD$ over $k$,
and choose non-constant, separable rational functions $f:\CC\to\P^1$
and $g:\DD\to\P^1$.  For simplicity, we assume that the zeroes and
poles of $f$ and $g$ are $k$-rational.  (This is not essential for any
of the results below, but it simplifies the formulas.)  Write the
divisors of $f$ and $g$ as
$$\dvsr(f)=\sum_{i=1}^k a_iP_i-\sum_{i'=1}^{k'} a'_{i'}P'_{i'}
\quad\text{and}\quad \dvsr(g)=\sum_{j=1}^\ell
b_jQ_j-\sum_{j'=1}^{\ell'} b'_{j'}Q'_{j'}$$ 
with $a_i,a'_{i'},b_i,b'_{j'}$ positive integers and $P_i$, $P'_{i'}$,
$Q_j$, and $Q'_{j'}$ distinct $k$-rational points.  (We hope the use
of $k$ for the ground field and for the number of zeroes of $f$ will
not cause confusion.)  Let
$$m=\sum_{i=1}^ka_i=\sum_{i'=1}^{k'} a'_{i'}
\quad\text{and}\quad
n=\sum_{j=1}^\ell b_j=\sum_{j'=1}^{\ell'} b'_{j'}.$$

We make two standing assumptions on the multiplicities $a_i$,
$a'_{i'}$, $b_j$, and $b'_{j'}$:
\begin{multline}\label{eq:mult-hyp}
\gcd(a_1,\dots,a_k,a'_1,\dots,a'_{k'},b_1,\dots,
b_\ell,b'_1,\dots,b'_{\ell'})=1\\
\text{for all $i$, $j$, }char(k)\nodiv\gcd(a_i,b_j)
\text{ and for all $i'$, $j'$, }char(k)\nodiv\gcd(a'_{i',}b'_{j'}).
\end{multline}
Here $char(k)$ is the characteristic of the ground field $k$.

Define a rational map $\pi_0:\CC\times_k\DD\ratto\P^1_k$ by the
formula
$$\pi_0(x,y)=[f(x):g(y)]$$
or, in terms of the standard coordinate $t$ on $P^1_k$, $t=f(x)/g(y)$.
It is clear that $\pi_0$ is defined away from the ``base points''
$(P_i,Q_j)$ ($1\le i\le k$, $1\le j\le\ell$) and $(P'_{i'},Q'_{j'})$
($1\le i'\le k'$, $1\le j'\le\ell'$).  Let $U$ be $\CC\times_k\DD$
with the base points removed, so that $\pi_0$ defines a morphism
$\pi_0|_U:U\to\P^1_k$.

It is useful to think of $\pi_0$ as giving a pencil of divisors on
$\CC\times_k\DD$: for each closed point $t\in\P^1_k$, the Zariski
closure in $\CC\times_k\DD$ of $(\pi_0|_U)^{-1}(t)$ is a divisor on
$\CC\times_k\DD$ passing through the base points which we denote as
$\overline{\pi_0^{-1}(t)}$.  As $t$ varies we get a pencil of divisors
with base locus exactly $\{(P_i,Q_j),(P'_{i'},Q'_{j'})\}$.  A key
feature of this pencil is that it fibers over $0$ and $\infty$ are
unions of ``horizontal'' and ``vertical'' divisors, i.e.,
$$\overline{\pi_0^{-1}(0)}=
\left(\bigcup_{i=1}^k \{a_iP_i\}\times\DD\right)
\cup 
\left(\bigcup_{j'=1}^{\ell'} \CC\times \{b'_{j'}Q'_{j'}\}\right)$$
and
$$\overline{\pi_0^{-1}(\infty)}=
\left(\bigcup_{i'=1}^{k'} \{a'_{i'}P_{i'}\}\times\DD\right)
\cup
\left(\bigcup_{j=1}^{\ell} \CC\times \{b_jQ_{j}\}\right).$$

In the proof of the proposition below, we will construct a specific
blow-up $\XX_1$ of $\CC\times_k\DD$ such that the composed rational
map $\XX_1\to\CC\times_k\DD\ratto\P^1$ is a generically smooth
morphism $\pi_1:\XX_1\to\P^1$.  (There are many such $\XX_1$, but the
statement of the theorem below is independent of this choice.)  Let
$X_1\to\spec K$ be the generic fiber of $\pi_1$ so that $X_1$ is a
smooth curve over $K$.

For each positive integer $d$ prime to the characteristic of $k$, let
$r_d:\P^1_k\to\P^1_k$ be the morphism with $r_d^*(t)=u^d$ (i.e.,
$u\mapsto t=u^d$ in the standard affine coordinates $t$ and $u$) and
form the base change
\begin{equation*}
\xymatrix{\mathllap{\SS_d:=}
\P^1_k\times_{\P^1_k}\XX_1\ar[r]\ar[d]&\XX_1\ar[d]^{\pi_1}\\
\P^1_k\ar[r]_{r_d}&\P^1_k}.
\end{equation*}
The fiber product $\SS_d$ will usually not be smooth over $k$ (or even
normal) because both $r_d$ and $\pi_1$ have critical points over $0$
and $\infty$.  Let $\XX_d\to \SS_d$ be a blow-up of the normalization
of the fiber product so that $\XX_d$ is smooth over $k$.  (Again,
there are many such, but the statement of the theorem is independent
of the choice.)  Let $\pi_d:\XX_d\to\P^1_k$ be the composition of
$\XX_d\to\SS_d$ with the projection onto $\P^1_k$.  (The domain
$\P^1_k$ for $\pi_d$ is thus the $\P^1_k$ in the lower left of the
displayed diagram.)  Let $K_d$ be the function field of this $\P^1_k$
so that $K_d=k(t^{1/d})\cong k(u)$ is an extension of $K$ of degree
$d$.  Finally, let $X_d\to\spec K_d$ be the generic fiber of $\pi_d$.
Then we have $X_d\cong X_1\times_{\spec K}\spec K_d$.

With this notation, we can state the main results of Berger's thesis
\cite{Berger08}*{2.2 and 3.1}:

\begin{thm}\label{thm:Berger}
  \textup{(Berger)} Choose data $k$, $\CC$, $\DD$, $f$, and $g$ as
  above, subject to the hypotheses~\ref{eq:mult-hyp}.  For each $d$
  prime to the characteristic of $k$, construct the fibered surfaces
  $\pi_d:\XX_d\to\P^1$ and the curves $X_d/K_d$ as above.
\begin{enumerate}
\item $X=X_1$ is a
  smooth, proper curve of genus
%\begin{multline*}
  $$\qquad\qquad g=mg_{\DD}+ng_{\CC}+(m -1)(n -1)%\\
  -\sum_{i,j}\delta(a_i,b_j)
  -\sum_{i',j'}\delta(a'_{i'},b'_{j'})
$$%\end{multline*}
where $g_{\CC}$ and $g_{\DD}$ are the genera of $\CC$ and $\DD$
respectively, $m=\sum_{i=1}^ka_i$, $n=\sum_{j=1}^\ell b_j$, and
$\delta(a,b)=(ab-a-b+\gcd(a,b))/2$.
\item $X_d$ is irreducible and remains irreducible over
  $\overline{k}K_d$.
\item $\XX_d$ is dominated by a product of curves.
\end{enumerate}
\end{thm}

\begin{rems} \mbox{}
\begin{enumerate}
\item Note that the generic fiber of $\pi_0$ is an open subset of the
  curve called $Y$ in the introduction and so the curve $X$ of the
  theorem is a smooth proper model of $Y$.
\item We will sketch the proof of part (1) of the proposition,
  essentially along the lines of Berger's thesis, in the rest of this
  section.
\item Berger proved that $\XX_d$ is DPC by using a fundamental group
  argument, generalizing \cite{Schoen96}*{6.6}.  In the next section,
  we will give a different and more explicit proof which will be the
  basis for a general rank formula.
\item Berger proved a stronger irreducibility result, namely that $X$
  is absolutely irreducible.  The result stated above is sufficient
  for our purposes and follows easily from the same ideas used to
  prove the domination by a product of curves.
\end{enumerate}
\end{rems}

\subsection{Proof of Thm.~\ref{thm:Berger}, part (1)}\label{ss:BergerBlowUps}
What we have to prove is that there is a blow-up of $\CC\times_k\DD$
such that $\pi_0$ becomes a generically smooth morphism, and that its
generic fiber has the properties announced in part (1).

Recall that $\pi_0|_U:U\to\P^1$ is the maximal open where $\pi_0$ is a
morphism.  An easy application of the Jacobian criterion shows that
$\pi_0|_U$ has only finitely many critical values and so all but
finitely many fibers of $\pi_0|_U$ are smooth (open) curves.  But
their closures $\overline{\pi_0^{-1}(t)}$ may be quite singular at the
base points.

Let us focus attention near a base point $(P,Q)=(P_i,Q_j)$.  (The
argument near a base point $(P'_{i'},Q'_{j'})$ is essentially
identical and will be omitted.)  Choose uniformizing parameters $x$
and $y$ on $\CC$ and $\DD$ at $P$ and $Q$ respectively.  Then in a
Zariski open neighborhood of $(P,Q)$, the map $\pi_0$ is defined away
from $(P,Q)$ and given by
$$(\alpha,\beta)\mapsto[x(\alpha)^a:y(\beta)^bu(\alpha,\beta)]$$
where $(a,b)=(a_i,b_j)$ and where $u$ is a unit in the local ring at
$(P,Q)$.  We say more succinctly that $\pi_0$ is given by
$[x^a:y^bu]$.  Now blowing up $(P,Q)$ and passing to suitable
coordinates, $\tilde\pi_0$, the composition of the blow-up with
$\pi_0$, is given by $[\tilde x^{a-b}:\tilde y^b\tilde u]$ if $a\ge b$
or by $[\tilde x^a:\tilde y^{b-a}\tilde u]$ if $a\le b$.  Thus if
$a=b$, $\tilde\pi_0$ is a morphism in a neighborhood of the inverse
image of $(P,Q)$. If $a\neq b$, there is exactly one point of
indeterminacy upstairs.  We relabel $\tilde\pi_0$ as $\pi_0$ and
iterate.  After finitely many blow-ups, $\pi_0$ is given by $[x^c,u]$
where $c=\gcd(a,b)$, and is defined everywhere on an open neighborhood
of the inverse image of $(P,Q)$.  We note for future use that the
exceptional fibers of the blow-ups at stages where $a\neq b$ map to
$0$ or $\infty\in\P^1_k$ whereas the exceptional divisor at the last
stage (when $a=b=c$) maps $c$-to-1 to $\P^1_k$.

For a pair of positive integers $(a,b)$, let $\gamma(a,b)$ be the
number of steps to proceed from $(a,b)$ to $(\gcd(a,b),0)$ by
subtracting the smaller of the pair from the larger at each step.
(So, for example, $\gamma(1,1)=1$, $\gamma(a,a)=1$, $\gamma(a,1)=a$
and $\gamma(2,3)=3$.)  The discussion above shows that if we blow up
$\CC\times_k\DD$ at the base points $(P_i,Q_j)$ $\gamma(a_i,b_j)$
times and at $(P'_{i'},Q'_{j'})$ $\gamma(a'_{i'},b'_{j'})$ times we
arrive at a surface $\XX_1\to\CC\times_k\DD$ such that the composed
rational map $\XX_1\to\CC\times_k\DD\ratto\P^1$ is a morphism
$\pi_1:\XX_1\to\P^1$.  I claim that $\pi_1$ is generically smooth with
fibers of genus as stated in the Proposition.

We have already seen that most fibers of $\pi_1$ are smooth away from
the base points.  Near a base point, $\FF=\overline{\pi_0^{-1}(t)}$ is
given locally by $tx^a-y^bu=0$.  If $a\neq b$ then the tangent cone is
a single line and the proper transform of $\FF$ has a single point
over the base point, at $x=y=0$ in the coordinates used above.
Continuing in this fashion, there is always a unique point over
$(P,Q)$ in the proper transform of $\FF$, except at the penultimate
step, where the proper transform of $\FF$ is given locally by
$tx^c-y^cu=0$ and its tangent cone consists of $c$ distinct lines.
(Here we use the hypothesis~\ref{eq:mult-hyp} which implies that $c$
is prime to the characteristic of $k$.)  In the final blow-up, the
proper transform has $c$ points over the base point and $\pi$ is
smooth in a neighborhood of these points, as long as we avoid $t=0$
and $t=\infty$.

This proves that all but finitely many fibers of $\pi_1$ are smooth
proper curves.  To find their genus, note that the fibers
$\overline{\pi_0^{-1}(t)}$ are curves of bidegree $(m,n)$ on
$\CC\times_k\DD$ and thus have arithmetic genus $mg_{\DD}+ng_{\CC}+(m
-1)(n -1)$.  A blow-up when the local equation is $tx^a-y^bu=0$
decreases the arithmetic genus by $e(e-1)/2$ where $e=\min(a,b)$.  The
value of the genus of the fibers of $\pi_1$ then follows by a simple
calculation.  (See \cite{Berger08}*{\S3.7 and \S3.8} for more
details.)

This completes the proof of part (1) of Theorem~\ref{thm:Berger}.
\qed

It is worth noting that the morphism $\pi_1:\XX_1\to\P^1_k$
constructed above may not be relatively minimal.

\section{Explicit domination by a product of curves}\label{s:eDPC}

\subsection{Kummer covers of $\CC$ and $\DD$}\label{ss:CC_d}
For each $d$, let $\CC_d$ be the smooth, proper, (possibly reducible)
covering of $\CC$ of degree $d$ defined by the equation $z^d=f(x)$ and
let $\DD_d$ be defined similarly by $w^d=g(y)$.  Let $e=e_{d,f}$ be
the largest divisor of $d$ such that $f$ is an $e$-th power in
$k(\CC)^\times$, so that $\CC_d$ is irreducible over $k$ if and only
if $e_{d,f}=1$.  Note that $e$ is a divisor of
$\gcd(d,a_1,\dots,a_k,a'_1,\dots,a'_{k'})$, with equality when
$k=\kbar$ and $\CC=\P^1$.  A similar discussion applies to $g$ and
$\DD_d$.  The hypothesis~\ref{eq:mult-hyp} implies that
$\gcd(e_{d,f},e_{d,g})=1$.

Let $\CC_d^o$ be $\CC_d$ with the (inverse images of the) zeroes and
poles of $f$ removed, and let $\DD_d^o$ be $\DD_d$ with the (inverse
images of the) zeroes and poles of $g$ removed.  Then
$\CC_d^o\to\CC_1^o$ and $\DD_d^o\to\DD_1^o$ are \'etale torsors for
$\mu_d$.  In particular, there is a free action of the \'etale group
scheme $\mu_d$ over $k$ on $\CC_d^o$ and $\DD_d^o$ and the quotients
are $\CC_1^o$ and $\DD_1^o$.

\subsection{Analysis of a fiber product}
Let $\SS=\CC\times\DD$ and let $D\subset\SS$ be the reduced divisor
supported on $\CC\times\DD\setminus\CC^o\times\DD^o$.  Let
$\phi_1:\XX_1\to\SS$ be the blow-up described in the proof of
Theorem~\ref{thm:Berger}~(1).  For $d>1$, let $\SS_d$ be the fiber
product $\P^1_k\times_{\P^1_k}\XX_1$ defined in the previous section
and let $\phi_d:\XX_d\to\SS$ be the composition
$\XX_d\to\SS_d\to\XX_1\to\SS$.  The following diagram may be helpful
in organizing the definitions:
\begin{equation}\label{eq:fiber-product}
\xymatrix{
\XX_d\ar[r]\ar[rdd]_{\pi_d}&\SS_d\ar[r]\ar[dd]
     &\XX_1\ar[d]_{\phi_1}\ar@/^4pc/[dd]^{\pi_1}\\
&&D\subset\CC\times_k\DD=\SS\ar@{-->}[d]\\
&\P^1_k\ar[r]_{r_d}&\P^1_k}
\end{equation}

For $d\ge1$, we define $\XX_d^o$ to be $\XX_d$ minus the divisor
$\phi_d^{-1}(D)$.  It is clear that $\phi_1$ induces an isomorphism
$\XX_1^o\cong\CC^o\times\DD^o$ and that $\XX_d^o\cong\SS_d^o$ is
isomorphic to the fiber product of
$\XX_1^o\to\P^1\setminus\{0,\infty\}$ with the \'etale morphism
$\P^1\setminus\{0,\infty\}\to\P^1\setminus\{0,\infty\}$, $t\mapsto
t^d$.

\begin{prop}\label{prop:explicitDPC}
There is a canonical isomorphism
$$\XX_d^o\cong\left(\CC_d^o\times\DD_d^o\right)/\mu_d$$
where $\mu_d$ acts on $\CC_d^o\times\DD_d^o$ diagonally.
\end{prop}

\begin{proof}
  By definition, $\XX_1^o$ is the open subset of $\CC\times_k\DD$
  where $f(x)$ and $g(y)$ are $\neq0,\infty$.  By the discussion above,
  $\XX_d^o$ is the closed subset of
$$\XX_1^o\times_k\A_k^1=\CC^o\times_k\DD^o\times_k\A^1_k$$
(with coordinates $(x,y,t)$) where $f(x)=t^dg(y)$.

On the other hand, $\CC_d^o\times_k\DD_d^o$ is isomorphic to the
closed subset of
$$\CC^0\times_k\A^1_k\times_k\DD^o\times_k\A^1_k$$
(with coordinates $(x,z,y,w)$) where $f(x)=z^d$ and $g(y)=w^d$.  It is
then clear that the morphism $(x,z,y,w)\mapsto(x,y,z/w)$ presents
$\CC_d^o\times\DD_d^o$ as a $\mu_d$-torsor over $\XX_d^o$.
\end{proof}

\subsection{Proof of Thm.~\ref{thm:Berger}, parts (2) and (3)}
\label{ss:end-of-proof}
Since $\XX_d^o$ is birational to $\XX_d$,
Proposition~\ref{prop:explicitDPC} shows that $\XX_d$ is dominated by
a product of curves, i.e., we have part (3) of
Theorem~\ref{thm:Berger}.

For part (2), we may assume that $k=\kbar$.  In this case, $e_{d,f}$
and $e_{d,g}$ (defined in the first paragraph of
Subsection~\ref{ss:CC_d}) are the numbers of irreducible components of
$\CC_d$ and $\DD_d$ respectively.  The components of $\CC_d$ are a
torsor for $\mu_{e_{d,f}}$ those of $\DD_d$ are a torsor
for $\mu_{e_{d,g}}$.  The action of
$\mu_d$ on the components is via the quotients
$\mu_d\to\mu_{e_{d,f}}$ and $\mu_d\to\mu_{e_{d,g}}$.  By
hypothesis~\ref{eq:mult-hyp}, $\gcd(e_{d,f},e_{d,g})=1$.  This implies
that under the diagonal action, $\mu_d$ acts
transitively on the set of components of $\CC_d\times_k\DD_d$ and so
$(\CC_d\times_k\DD_d)/\mu_d$ is irreducible.  It follows
that $\XX_d$ is geometrically irreducible and therefore that $X_d$ is
irreducible over $\overline{k}K_d$.

This completes the proof of Theorem~\ref{thm:Berger}.  \qed

This proof suggests another way to construct $\XX_d$, namely starting
from $\CC_d\times_k\DD_d$.  In the next section we will carry out this
construction and use it to give a formula for the rank of
$MW(J_{X_d})$.

\subsection{The $K/k$ trace of $J_d$}
We end this section with an analysis of $B_d$, the $K/k$-trace of
$J_d=Jac(X_d)$ under the simplifying assumption that $k$ is
algebraically closed.  Fix $d$ and recall $e_{d,f}$ and $e_{d,g}$ from
Subsection~\ref{ss:CC_d}.  Fix also
roots $f'=f^{1/e_{d,f}}$ and $g'=g^{1/e_{d,g}}$ and write $\CC'_{d'}$ and
$\DD'_{d'}$ for the covers of $\CC$ and $\DD$ defined by the equations
$z^{d'}=f'$ and $w^{d'}=g'$.

\begin{prop}\label{prop:B}
  Assume that $k$ is algebraically closed. With notation as above, 
  the $K/k$-trace $B_d$ of $J_d=Jac(X_d)$ is canonically
  isomorphic to $J_{\CC'_{e_{d,g}}}\times_k J_{\DD'_{e_{d,f}}}$.
\end{prop}

(The positions of $e_{d,f}$ and $e_{d,g}$ in this formula are not
typos---see the proof below.)

\begin{proof}
  By Proposition~\ref{prop:ST}, $B_d$ is isomorphic to the Picard
  variety of $\XX_d$ and since the Picard variety is a birational
  invariant, we need to compute
  $\picvar\left((\CC_d\times_k\DD_d)/\mu_d\right)$.  Applying
  \cite{Shioda99}*{Thm~2} twice, we have that for irreducible curves
  $\CC$ and $\DD$ over $k$, $\picvar(\CC\times_k\DD)\cong
  J_\CC\times_k J_\DD$.

  Now the components of $\CC_d$ (resp.~$\DD_d$) are indexed by
  $\mu_{e_{d,f}}$ (resp.~$\mu_{e_{d,g}}$) and each one is isomorphic to
  $\CC'_{d/e_{d,f}}$ (resp.~$\DD'_{d/e_{d,g}}$).  The action of $\mu_d$
  on the components $\CC_d\times_k\DD_d$ is transitive, with
  stabilizer $\mu_{d/(e_{d,f}e_{d,g})}$.  Thus
$$\picvar\left((\CC_d\times_k\DD_d)/\mu_d\right)
\cong\picvar\left((\CC'_{d/e_{d,f}}\times_k\DD'_{d/e_{d,g}})
  /\mu_{d/(e_{d,f}e_{d,g})}\right).$$ 
Applying the formula above for products of curves and passing to the
$\mu_{d/(e_{d,f}e_{d,g})}$ invariants, we have
$$B_d\cong J_{\CC'_{e_{d,g}}}\times_k J_{\DD'_{e_{d,f}}}.$$
This completes the proof of the proposition.
\end{proof}

A nice example occurs in \cite{Berger08}*{4.3(4)} where
$\CC=\DD=\P^1_k$, $f(x)=x^2$ and $g(y)$ is a quadratic rational
function with distinct zeroes and poles so that $e_{d,f}=\gcd(2,d)$,
$e_{d,g}=1$.  It follows that when $d$ is odd, $B_d=0$ and when $d$ is
even, $B_d$ is the elliptic curve $w^2=g(y)$.  This is confirmed in
\cite{Berger08} by explicit computation.

\begin{rem}
Let $B'_d$ be the image of the $K/k$-trace $\tau:B_d\times_kK\to J_d$
and define an abelian variety $A_d$ so that
$$0\to B'_d\to J_d\to A_d\to0$$
is exact.  Since $\tau$ is purely inseparable \cite{Conrad06}*{6.12},
we have $B_d'(K)=B_d(K)$ and since $K=k(\P^1_k)$, we have $B_d(K)=B_d(k)$.
Thus the Mordell-Weil group $J_d(K)/\tau B_d(k)$ is a subgroup of the
more traditional Mordell-Weil group $A_d(K)$.
\end{rem}

\section{A rank formula}\label{s:Formula}

\subsection{Numerical invariants}
Throughout this section, we assume that $k$ is algebraically closed.
Fix data $\CC$, $\DD$, $f$, and $g$ as in Berger's construction
(Subsection~\ref{ss:SpecialPencils}), subject to the
hypothesis~\ref{eq:mult-hyp}.  For each positive integer $d$ not
divisible by the characteristic of $k$, form covers $\CC_d\to\CC$ and
$\DD_d\to\DD$ defined by the equations $z^d=f$ and $w^d=g$.  Recall
that (since $k=\kbar$) $e_{d,f}$ and $e_{d,g}$ are the numbers of
irreducible components of $\CC_d$ and $\DD_d$ respectively.

For each $i$ (resp.~$i'$, $j$, $j'$), let $r_{d,i}$
(resp.~$r'_{d,i'}$, $s_{d,j}$, $s'_{d,j'}$) be the number of closed
points of $\CC_d$ over $P_i$ (resp.~of $\CC_d$ over $P'_{i'}$, of
$\DD_d$ over $Q_j$, of $\DD_d$ over $Q'_{j'}$).  Also, for each pair
$(i,j)$ (resp.~$(i',j')$) let $t_{d,i,j}$ (resp.~$t'_{d,i',j'}$) be the
number of closed points of $(\CC_d\times_k\DD_d)/\mu_d$ (quotient by
the diagonal action) over $(P_i,Q_j)$ (resp.~over
$(P'_{i'},Q'_{j'})$).  Since $k$ is algebraically closed, we have
equalities 
\begin{align*}
r_{d,i}&=\gcd(a_i,d)&r'_{d,i'}&=\gcd(a'_{i'},d)\\
s_{d,j}&=\gcd(b_j,d)&s'_{d,j'}&=\gcd(b'_{j'},d)\\
t_{d,i,j}&=\gcd(a_i,b_j,d)&t'_{d,i',j'}&=\gcd(a'_{i'},b'_{j'},d).
\end{align*}

On $\CC_d$, $z$ has a zero of order $a_i/r_{d,i}$ at each point
over $P_i$ and a pole of order $a'_{i'}/r'_{d,i'}$) at each
point over $P'_{i'}$.  Similarly for $w$ and $\DD_d\to\DD$. Let
\begin{align*}
&k_d=\sum_{i=1}^kr_{d,i}\qquad
k'_d=\sum_{i'=1}^{k'}r'_{d,i'}\\
&\ell_d=\sum_{j=1}^\ell s_{d,j}\qquad
\ell'_d=\sum_{j'=1}^{\ell'}s'_{d,j'}.
\end{align*}

\subsection{A nice model for $X_d$}
Now consider the rational map $\CC_d\times\DD_d\ratto\P^1_k$ defined
by $z/w$.  This is defined away from $k_d\ell_d+k'_d\ell'_d$ base
points.  Blowing up each one as in Subsection~\ref{ss:BergerBlowUps}
several times (more precisely, $\gamma(a_i/r_{d,i},b_j/s_{d,j})$ times
at points over $(P_i,Q_j)$ and $\gamma(a'_{i'}/r'_{d,i'},b'_{j'}/s'_{d,j'})$ 
times at points over $(P'_{i'},Q'_{j'})$), we arrive at a surface
$\widetilde{\CC_d\times_k\DD_d}$ equipped with a generically
smooth morphism to $\P^1_k$ which we again denote $z/w$.  Over
each base point, all but one component of the exceptional locus
maps to $0$ or $\infty$ in $\P^1_k$ and the remaining component
(the exceptional divisor of the last blow-up at that base point)
maps finitely to $\P^1_k$.  Write
$$\gamma_{d,i,j}=
\gamma(\frac{a_i}{r_{d,i}},\frac{b_j}{s_{d,j}})$$
($i=1,\dots,k$, $j=1,\dots,\ell$) and 
$$\gamma'_{d,i',j'}=
\gamma(\frac{a'_{i'}}{r'_{d,i'}},\frac{b'_{j'}}{s'_{d,j'}})$$
($i'=1,\dots,k'$, $j'=1,\dots,\ell'$) for the number of blow-ups
needed at each base point.

The diagonal action of $\mu_d$ on $\CC_d\times_k\DD_d$ lifts
canonically to the blow-up $\widetilde{\CC_d\times_k\DD_d}$.  The
locus where the action has non-trivial stabilizers has divisorial
components and sometimes has isolated fixed points (depending on the
order of zeroes or poles of $z$ and $w$ at the base points).  The
isolated fixed points all map via $z/w$ to $0$ or $\infty$ in
$\P^1_k$.  Let $\YY_d$ be the quotient
$(\widetilde{\CC_d\times_k\DD_d})/\mu_d$.  The morphism
$z/w:\widetilde{\CC_d\times_k\DD_d}\to\P^1_k$ descends to a morphism
which is generically smooth and which we again denote
$z/w:\YY_d\to\P^1$.  If the action of $\mu_d$ on
$\widetilde{\CC_d\times_k\DD_d}$ has any isolated fixed points, then
$\YY_d$ has cyclic quotient singularities as studied by Hirzebruch and
Jung.  The resolution of these singularities is well known: the
exceptional locus is a chain of $\P^1$s whose length and
self-intersections are given by a Hirzebruch-Jung continued fraction
(see \cite{BHPV}*{III.5}).  Let $\XX_d\to\YY_d$ be the canonical
resolution of the quotient singularities and let $HJ_d$ be the number
of irreducible components of the exceptional locus.
%use Beauville VIII.5 to justify that these components
%are independent in NS so the rank contribution is also HJ_d:
%the components obviously generate the extra NS and there are no
% relations because there are none even numerically
The only point about the singularities that will be relevant for us is
that every component of the exceptional locus lies over $0$ or
$\infty\in \P^1_k$.  Write $\pi_d$ for the composition of
$\XX_d\to\YY_d$ with $z/w$.  Then $\XX_d$ is smooth and proper over
$k$, $\pi_d:\XX_d\to\P^1_k$ is generically smooth, and $\XX_d$ is
birational to the fiber product $\SS_d$ of
Section~\ref{ss:SpecialPencils}.  So $\XX_d$ is a nice model of $X_d$
and we can use it to compute the rank of $MW(J_d)$.

\subsection{More numerical invariants}
To state our main result, we define two more numerical invariants.
For each closed point $v$ of $\P^1_k$, let $f_{d,v}$ be the number of
irreducible components in the fiber of $\pi_d$ over $v$.  Define
$$c_1(d)=\sum_{v\neq0,\infty}(f_{d,v}-1)=d\sum_{v\neq0,\infty}(f_{1,v}-1).$$
Here the first sum is over places of the $\P^1$ in the lower left
corner of diagram~\ref{eq:fiber-product} and the second is over places
of the $\P^1$ in the lower right of the same diagram.  The second
equality uses that $k$ is algebraically closed.

Define a second function $c_2(d)$ by
\begin{multline}
c_2(d)=\sum_{i=1,\dots,k,j=1,\dots,\ell}t_{d,i,j}
-\sum_{i=1,\dots,k}\gcd(a_i,e_{d,g})-\sum_{j=1,\dots,\ell}\gcd(b_j,e_{d,f})
+1\\
+\sum_{i'=1,\dots,k',j'=1,\dots,\ell'}t'_{d,i',j'}
-\sum_{i'=1,\dots,k'}\gcd(a'_{i'},e_{d,g})-\sum_{j'=1,\dots,\ell'}\gcd(b'_{j'},e_{d,f})
+1.
\end{multline}
For $d$ prime to the multiplicities $a_i$, $a'_{i'}$, $b_j$,
$b'_{j'}$, we have
$$c_2(d)=c_2(1)=(k-1)(\ell-1)+(k'-1)(\ell'-1).$$ 
For all $d$, we have $c_2(d)\ge c_2(1)$ and $c_2(d)$ is periodic and
so bounded.

We note that $c_2$ is of a combinatorial nature---it depends only on
the multiplicities $a_i$, $a'_{i'}$, $b_j$, and $b'_{j'}$---whereas
$c_1$ in general depends also on the positions of the points $P_i$,
$P'_{i'}$, $Q_j$, and $Q'_{j'}$.

Before stating the theorem, we recall that $e_{d,f}$ is the largest
divisor of $d$ such that $f$ is an $e_{d,f}$-th power in
$k(\CC)^\times$ and similarly for $e_{d,g}$.  Choose roots
$f'=f^{1/e_{d,f}}$ and $g'=g^{1/e_{d,g}}$ and let $\CC'_{d/e_{d,f}}$
and $\DD'_{d/e_{d,g}}$ be the covers of $\CC$ and $\DD$ defined by
$z^{d/e_{d,f}}=f'$ and $w^{d/e_{d,g}}=g'$ respectively.  Their
Jacobians will be denoted $J_{\CC'_{d/e_{d,f}}}$ and
$J_{\DD'_{d/e_{d,g}}}$ respectively.

\begin{thm}\label{thm:main}
  Assume that $k$ is algebraically closed.  Choose data $\CC$, $\DD$,
  $f$, and $g$ as above, subject to the hypotheses~\ref{eq:mult-hyp}.
  Let $X$ be the smooth proper model of
$$\{f-tg=0\}\subset\CC\times_k\DD\times_k\spec K$$  
over $K=k(t)$ constructed in Section~\ref{s:Berger}.  For each $d$
prime to the characteristic of $k$ let $K_d=k(t^{1/d})$, and let $J_d$
be the Jacobian of $X$ over $K_d$.  Then with notation as above,
$$\rk MW(J_d)=
\rk\Hom_{k-av}(J_{\CC'_{d/e_{d,f}}},J_{\DD'_{d/e_{d,g}}})^{\mu_{d/e_{d,f}e_{d,g}}}
-c_1(d)+c_2(d).$$ 
Here $\Hom_{k-av}$ denotes homomorphisms of abelian varieties over $k$
and the exponent $\mu_{d/e_{d,f}e_{d,g}}$ signifies those
homomorphisms which commute with the action of
$\mu_{d/e_{d,f}e_{d,g}}$.
\end{thm}

In the case where $e_{d,f}=e_{d,g}=1$, so that $\CC_d$ and $\DD_d$ are
irreducible, the rank formula simplifies to
$$\rk MW(J_d)=
\rk\Hom_{k-av}(J_{\CC_{d}},J_{\DD_{d}})^{\mu_{d}}
-c_1(d)+c_2(d).$$

\begin{proof}
  The proof consists of computing the rank of the N\'eron-Severi group
  of $\XX_d$ using the construction of $\XX_d$ via
  $\CC_d\times_k\DD_d$ and then applying the Shioda-Tate formula.

  First note that $\CC_d\times_k\DD_d$ is isomorphic to the disjoint
  union of $e_{d,f}e_{d,g}$ copies of
  $\CC'_{d/e_{d,f}}\times_k\DD'_{d/e_{d,g}}$ so the N\'eron-Severi
  group of $\CC_d\times_k\DD_d$ is isomorphic to
$$\left(\Hom_{k-av}(J_{\CC'_{d/e_{d,f}}},J_{\DD'_{d/e_{d,g}}})\oplus\Z^2\right)^{e_{d,f}e_{d,g}}.$$

Each blow-up in passing from $\CC_d\times_k\DD_d$ to
$\widetilde{\CC_d\times_k\DD_d}$ contributes a factor of $\Z$ to
N\'eron-Severi.  Grouping these by the base point $(P_i,Q_j)$ or
$(P'_{i'},Q'_{j'})$ over which they lie, we have
\begin{align*}
NS(\widetilde{\CC_d\times_k\DD_d})&\cong
\left(\Hom_{k-av}(J_{\CC'_{d/e_{d,f}}},J_{\DD'_{d/e_{d,g}}})\oplus\Z^2\right)^{e_{d,f}e_{d,g}}\\
&\quad\bigoplus_{i,j}\Z^{r_{d,i}s_{d,j}\gamma_{d,i,j}}\\
&\quad\bigoplus_{i',j'}\Z^{r'_{d,i'}s'_{d,j'}\gamma'_{d,i',j'}}.
\end{align*}

Now consider the action of $\mu_d$ on
$NS(\widetilde{\CC_d\times_k\DD_d})$.  The group $\mu_d$ acts
transitively on the components of $\CC_d\times_k\DD_d$ and the
stabilizer of a component is $\mu_{d/e_{d,f}e_{d,g}}$.  Thus the
invariants under $\mu_d$ of the first summand above is
$$\Hom_{k-av}(J_{\CC'_{d/e_{d,f}}},J_{\DD'_{d/e_{d,g}}})^{\mu_{d/e_{d,f}e_{d,g}}}
\oplus\Z^2$$ 
where we use an exponent to denote invariants.  For a fixed $(i,j)$,
the action of $\mu_d$ on the base points of $z/w$ on
$\CC_d\times_k\DD_d$ over $(P_i,Q_j)\in\CC\times_k\DD$ has $t_{d,i,j}$
orbits.  The stabilizer of one of these points acts trivially on the
part of the exceptional locus of
$\widetilde{\CC_d\times_k\DD_d}\to\CC_d\times_k\DD_d$ over that point.
Thus the contribution to the rank of
$NS(\widetilde{\CC_d\times_k\DD_d}/\mu_d)$ corresponding to
$(P_i,Q_j)$ is $t_{d,i,j}\gamma_{d,i,j}$.  A similar discussion
applies to the base points $(P'_{i'},Q'_{j'})$.  Thus the rank of the
N\'eron-Severi group of $\widetilde{\CC_d\times_k\DD_d}/\mu_d$ is
$$\rk\Hom_{k-av}(J_{\CC'_{d/e_{d,f}}},J_{\DD'_{d/e_{d,g}}})^{\mu_{d/e_{d,f}e_{d,g}}}
+2
+\sum_{i,j}t_{d,i,j}\gamma_{d,i,j}
+\sum_{i',j'}t'_{d,i',j'}\gamma'_{d,i',j'}.$$ 
The resolution of the Hirzebruch-Jung singularities contributes an
additional $HJ_d$ (by definition) to the rank of $NS(\XX_d)$, so we
have
\begin{align}\label{eq:NS(Xd)}
\rk NS(\XX_d)&=
\rk\Hom_{k-av}(J_{\CC'_{d/e_{d,f}}},J_{\DD'_{d/e_{d,g}}})^{\mu_{d/e_{d,f}e_{d,g}}}
+2\\
&\quad+\sum_{i,j}t_{i,j}\gamma_{d,i,j}
+\sum_{i',j'}t'_{d,i',j'}\gamma'_{d,i',j'}\notag\\
&\quad+HJ_d.\notag
\end{align}

Now consider the morphism $\pi_d:\XX_d\to\P^1_k$ and the Shioda-Tate
formula.  There are three sources of components of $\XX_d$ over 0 or
$\infty$ in $\P^1_k$: (i) the images in $\XX_d$ of components of
curves on $\CC_d\times_k\DD_d$ over the curves $\{P_i\}\times_k\DD$,
$\{P_{i'}\}\times_k\DD$, $\CC\times_k\{Q_j\}$, and $\CC\times_k\{Q_{j'}\}$
in $\CC\times_k\DD$; (ii) the images in $\XX_d$ of the exceptional
divisors of the blow-ups needed to pass from $\CC_d\times_k\DD_d$ to
$\widetilde{\CC_d\times_k\DD_d}$ (more precisely the exceptional
divisors for all but the last blow up at each base point---see the
third paragraph of Subsection~\ref{ss:end-of-proof}); and (iii) the
exceptional divisors of the blow-ups needed to resolve any
Hirzebruch-Jung singularities.
Taking into account the action of $\mu_d$
as in the previous paragraph we find that
\begin{align*}
f_{d,\infty}+f_{d,0}
&=\sum_{i=1,\dots,k}\gcd(a_i,e_{d,g})+\sum_{i'=1,\dots,k'}\gcd(a'_{i'},e_{d,g})\\
&+\sum_{j=1,\dots,\ell}\gcd(b_j,e_{d,f})+\sum_{j'=1,\dots,\ell'}\gcd(b'_{j'},e_{d,f})\\
&\quad+\sum_{ij}t_{d,i,j}(\gamma_{d,i,j}-1)+
\sum_{i'j'}t'_{d,i',j'}(\gamma'_{d,i',j'}-1)\\
&\quad+HJ_d.
\end{align*}
Combining this and equation~\ref{eq:NS(Xd)} with the Shioda-Tate
formula~\ref{prop:ST} yields the stated formula for $\rk MW(J_d)$.

This completes the proof of Theorem~\ref{thm:main}
\end{proof}

\begin{rems}\label{rem:data}\mbox{}
\begin{enumerate}
\item We used the hypothesis that $k$ is algebraically closed in two
  ways.  First, it is used in the Shioda-Tate formula.  In fact, this
  formula is known to hold in much greater generality; this will be
  discussed in detail in a forthcoming publication.  Second, without
  this hypothesis, rationality issues would come into the definitions
  of the numerical invariants.  This is more of an inconvenience than
  a fundamental issue.  Thus, with some adjustments, we have a rank
  formula without the assumption that $k$ is algebraically closed.
\item The reader will note that there is a great deal of flexibility
  in choosing the input data for Berger's construction.  It is useful
  to think of first choosing the topological data, namely the genera
  of $\CC$ and $\DD$ and the multiplicities in the divisors of $f$ and
  $g$, and then consider the (connected) family of choices of the
  isomorphism classes of the curves and the locations of the zeroes
  and poles.
\item There is an interesting difference in the nature of the three
  terms in the rank formula.  The term $c_2$ depends only on the
  topological data.  The term $c_1$ depends on the continuous data in
  an algebro-geometric way---the sets on which it is constant are
  constructible.  The $\rk\Hom(\cdots)$ term is more arithmetic in
  nature---it may jump up on a countable collection of subschemes.
  Thus we may expect that the rank may go up countably often.  The
  examples in the remainder of the paper illustrate two facets of this
  remark.  One involves ranks going up in all positive characteristics
  and the other involves ranks going up for countably many values of a
  complex parameter.
\end{enumerate}
\end{rems}

In the rest of the paper we discuss two examples which illustrate
certain aspects of Berger's construction and the rank formula.
Further applications of these ideas appear in
\cite{OcchipintiThesis} and papers based on it.

\section{First example}\label{s:Example1}

The first interesting examples occur already when $\CC$ and $\DD$ are
rational and $f$ and $g$ are quadratic.

\subsection{Data}
Let $\CC=\DD=\P^1_k$, and let
$$f(x)=x(x-1)\quad\text{and}\quad g(y)=1/f(1/y)=y^2/(1-y).$$
The hypotheses~\ref{eq:mult-hyp} on $f$ and $g$ are satisfied no
matter what the characteristic of $k$, and for all $d$,
$e_{d,f}=e_{d,g}=1$. 

\subsection{The surface $\XX_1$}
The surface $\XX_1$ is obtained by blowing up $\CC\times\DD$ twice
($\gamma(2,1)=2$) at each of the four base points $(x,y)=(0,0)$
$(1,0)$, $(\infty,1)$, and $(\infty,\infty)$.

Straightforward calculation shows that the morphism
$\pi_1:\XX_1\to\P^1_k$ is smooth away from $t=0$, $\infty$, and, if
$char(k)\neq2$, $t=1/16$.  The smooth fibers have genus 1.  The fiber
over 0 is a configuration of Kodaira type $I_4$, the fiber over
$\infty$ is a configuration of Kodaira type $I_1^*$, and the fiber
over $t=1/16$ is an irreducible rational curve with one node, i.e., a
configuration of Kodaira type $I_1$.

It follows immediately that $c_1(d)=c_2(d)=0$ for all $d$.

\subsection{The curve $X$}
The curve
$$Y=\{f-tg=0\}\subset\P^1_k\times_k\P^1_k\times\spec K
\cong\P^1_K\times_K\P^1_K$$ 
is smooth of bidegree $(2,2)$ and so has genus 1.  The change of
coordinates $x=-y/(x+t)$, $y=-x/t$ brings it into the Weierstrass form
$$X:\qquad y^2+xy+ty=x^3+tx^2.$$
The discriminant of this model is $\Delta=t^4(1-16t)$ and the
$j$-invariant is $j=(16t^2-16t+1)^3/\Delta$.  (Dick Gross points out
that over $\Z[1/2]$, $X$ is the universal elliptic curve over the
modular curve $X_1(4)\cong\P^1$.)

Since $X$ is elliptic, it is its own Jacobian and by
Proposition~\ref{prop:B} (or because $j(X)\not\in k$), the $K/k$ trace
of $X_d$ is 0 for all $d$.  Thus the rank formula~\ref{thm:main} says
that if $k$ is algebraically closed
$$\rk X_d(K_d)=\rk\Hom_{k-av}(J_{\CC_d},J_{\DD_d})^{\mu_d}.$$

\subsection{Homomorphisms and endomorphisms}
Note that $\CC_d\cong\DD_d$ via $(x,z)\mapsto(y=1/x,w=1/z)$.  This
isomorphism, call it $\sigma$, anti-commutes with the $\mu_d$ action in
the sense that $\sigma\circ[\zeta_d]=[\zeta_d^{-1}]\circ\sigma$.
Using it to identify $\CC_d$ with $\DD_d$, we have
$$\Hom_{k-av}(J_{\CC_d},J_{\DD_d})^{\mu_d}\cong
\en_{k-av}(J_{\CC_d})^{anti-\mu_d}$$
where the superscript denotes those endomorphisms anti-commuting
with the $\mu_d$ action.  Thus if $k$ is algebraically closed we have
$$\rk X_d(K_d)=\rk\en_{k-av}(J_{\CC_d})^{anti-\mu_d}.$$

We can make this much more explicit.  Note that since $\CC_d$ is
defined over the prime field, every endomorphism of $J_{\CC_d}$ is
defined over an algebraic extension of the prime field.  Thus there is
no loss in assuming that $k$ is the algebraic closure the prime field.

Let $\phi(e)$ be the cardinality of $(\Z/e\Z)^\times$ and let $o_e(q)$
be the order of $q$ in $(\Z/e\Z)^\times$.

\begin{thm}\label{thm:ex1}
  If $char(k)=0$, then $\rk X_d(K_d)=0$ for all $d$.  If
  $char(k)=p>0$, then $\rk X_d(K_d)$ is unbounded as $d$ varies.  More
  precisely, suppose $k=\Fq$ and $d=p^n+1$.  Then
$$\rk X_d(K_d)\ge\sum_{e|d, e>2}\phi(e)/o_e(q)\ge (p^n-1)/2n.$$
If $k$ contains the $d$-th roots of unity, $\rk
X_d(K_d)=p^n-1$ is $p$ is odd and $p^n$ if $p=2$.
\end{thm}

The proof of the theorem will occupy the rest of this section.  In the
following section, we will see that the points asserted to exist here
can be made quite explicit.

\subsection{The case $char(k)=0$}
We begin with the case $char(k)=0$ where we may assume that $k$ is
algebraically closed.  Note that $\CC_d$ has genus $(d-1)/2$ when $d$
is odd and $(d-2)/2$ when $d$ is even.  Let $\zeta_d\in k$ be a
primitive $d$-th root of unity and $[\zeta_d]$ the corresponding
endomorphism of $J_{\CC_d}$.  We consider the faithful action of
$\en(J_{\CC_d})$ on $H^0(J_{\CC_d},\Omega^1)=H^0(\CC_d,\Omega^1)$.  It
is easy to write down an explicit basis of the 1-forms and to see that
the eigenvalues of $[\zeta_d]$ on this space are $\zeta_d^i$
$i=1,...,g_{\CC_d}$.  In particular, the eigenvalues of $[\zeta_d]$
and those of $[\zeta_d^{-1}]$ are disjoint.  Therefore, no non-zero
endomorphism $\phi$ of $J_{\CC_d}$ can satisfy
$\phi\circ[\zeta_d]=[\zeta_d^{-1}]\circ\phi$.  Thus
$$\rk X_d(K_d)=\rk \en_{k-av}(J_{\CC_d})^{anti-\mu_d}=0$$
and this proves the first part of the theorem.

\subsection{Old and new}
Between here and Subsection~\ref{ss:non-alg-closed-k}, we assume that
$k$ is algebraically closed of characteristic $p>0$.  The
modifications needed to treat non-algebraically closed $k$ will be
given in that subsection.

The discussion in this subsection is not specific to the example
under consideration and makes sense in the general context of Berger's
construction (Section~\ref{s:Berger}). 

If $e$ is a divisor of $d$ then there is a natural map $\CC_d\to\CC_e$
($z\mapsto z^{d/e}$) which induces an injective homomorphism
$J_{\CC_e}\to J_{\CC_d}$.  We write $A_d^{\text{old}}\subset
J_{\CC_d}$ for the sum over all $e|d$, $e\neq d$ of the images of
these maps.  Choose a complement to $A_d^{\text{old}}$ in $J_{\CC_d}$
(defined only up to isogeny), denoted $A_d^{\text{new}}$, so that we
have an isogeny $J_{\CC_d}\to A_d^{\text{old}}\oplus
A_d^{\text{new}}$.  It is clear that the action of $\mu_d$ on
$J_{\CC_d}$ preserves $A_d^{\text{old}}$ and therefore induces a
homomorphism $\mu_d\to\aut(A_d^{\text{new}})\subset\en(A_d^{\text{new}})$.

\begin{lemma}\label{lemma:endos}
  Let $G$ be a cyclic group of order $d$ and let $R$ be the group ring
  $\Q[G]$.  Let $R\to\en(J_{\CC_d})\tensor\Q$ be a homomorphism which
  sends a generator of $G$ to the action of a primitive $d$-th root of
  unity on $J_{\CC_d}$.  Suppose that $A_d^{\text{new}}$ is non-zero.
  Then the image of the composed map
$$R\to\en(J_{\CC_d})\tensor\Q\to\en(A_d^{\text{new}})\tensor\Q$$
is isomorphic to the field $\Q(\mu_d)$.
\end{lemma}

\begin{proof}
  When $k$ has characteristic zero, the lemma follow easily from a
  consideration of the action of $\zeta_d$ on 1-forms, but we need a
  proof that works over any $k$.  For this, \'etale cohomology serves
  as a good replacement for the 1-forms.

  Choose a prime $\ell$ not equal to the characteristic of $k$.
  Consider the (faithful) action of
  $\en^0(J_{\CC_d})=\en(J_{\CC_d})\tensor\Q$ on
$$H^1(J_{\CC_d},\Ql)=
\bigoplus_{e|d}H^1(A_e^{\text{new}},\Ql).$$
It will suffice to determine the image of $R$ in
$\en(H^1(A_d^{\text{new}},\Ql))$. 

If $r$ is an element of $R$, we will write $H^1(\cdots)^{r=0}$ for the
kernel of $r$ on $H^1(\cdots)$.

Let $\Phi_e(x)$ be the cyclotomic polynomial of order $e$, so that
$x^d-1=\prod_{e|d}\Phi_e(x)$.  Since $R\cong\Q[x]/(x^d-1)$, the
Chinese remainder theorem implies that 
$$R\cong\prod_{e|d}\Q(\mu_e).$$

Now the quotient of $\CC_d$ by the group generated by a primitive
$e$-th root of unity $\zeta_d^{d/e}$ is $\CC_e$.  This implies that
$$H^1(J_{\CC_d},\Ql)^{[\zeta_d^{d/e}]-1=0}
=H^1(J_{\CC_e},\Ql)\subset
H^1(A_d^{\text{old}},\Ql).$$

This in turn implies that in the map
$R\to\en(H^1(A_d^{\text{new}},\Ql)$, all the factors $\Q(\mu_e)$ with
$e<d$ go to zero.  But if $A_d^{\text{new}}\neq0$, then
$R\to\en(H^1(A_d^{\text{new}},\Ql)$ is not the zero homomorphism and
this forces its restriction to the factor $\Q(\mu_d)$ to be non-zero.
Since $\Q(\mu_d)$ is a field, it maps isomorphically onto its image
and this proves our claim.
\end{proof}

\subsection{Dimensions}
In the context of the data chosen at the beginning of this section,
the Riemann-Hurwitz formula implies that $A_e^{\text{new}}=0$ for
$e=1$ and $e=2$, and that it has dimension $\phi(e)/2$ when $e>2$.
The lemma implies that the image of $R$ in $\en^0(J_{\CC_d})$ has
dimension $(d-1)/2$ when $d$ is odd and $(d-2)/2$ when $d$ is even.

\subsection{Special endomorphisms}
Now suppose that $d$ has the special form $d=p^n+1$.  Let
$\Fr_q:\CC_d\to\CC_d$ be the $q$-power Frobenius where $q=p^n$.  Then
it is immediate that $[\zeta_d]\circ\Fr_q=\Fr_q\circ[\zeta_d^{-1}]$,
i.e., $\Fr_q$ gives an element of $\en(J_{\CC_d})^{anti-\mu_d}$.  The
same calculation shows that
$\Fr_q\circ[\zeta_d^i]\in\en(J_{\CC_d})^{anti-\mu_d}$ for
$i=0,\dots,d-1$.  Since $\Fr_q$ is not a zero-divisor in
$\en(J_{\CC_d})$, the relations among the $\Fr_q\circ[\zeta_d^i]$ are
the same as the relations among the $[\zeta_d^i]$.
Lemma~\ref{lemma:endos} shows that the endomorphisms $[\zeta_d^i]$
$i=0,\dots,d-1$ generate a subalgebra of $\en(J_{\CC_d})$ of rank
$p^n-1$ when $p$ is odd and of rank $p^n$ when $p=2$.  Thus the
endomorphisms $\Fr_q\circ[\zeta_d^i]$ $i=0,\dots,d-1$ span a
sublattice of $\en(J_{\CC_d})^{anti-\mu_d}$ of the same rank.  This
shows that when $k$ has characteristic $p$ and $d=p^n+1$, then the
rank of $X_d(K_d)$ is bounded below by $p^n-1$ when $p$ is odd and by
$p^n$ when $p=2$.

There are several ways to see that these lower bounds are
equalities. One is to examine more closely the structure of
$\en^0(J_{\CC_d}^{\text{new}})$.  Another is to reduce to the case
where $k$ is the algebraic closure of a finite field and then note
that by the Grothendieck-Ogg-Shafarevitch formula, the degree of the
$L$-function of $X_d$ over $K_d$ is $p^n-1$ ($p$ odd) or $p^n$ ($p=2$)
and so we have an {\it a priori\/} upper bound on the rank.  We omit
the details.

\subsection{General $k$}\label{ss:non-alg-closed-k}
For $k$ not necessarily algebraically closed, we can compute the rank
of $X_d(K_d)$ by passing to $\kbar$ and then taking invariants
under $\gal(\kbar/k)$.  The calculations are straightforward and so
again we omit the details.

\begin{rems} \mbox{}
\begin{enumerate}
\item The example above, at least as regards analytic ranks, can also
  be treated by the methods of \cite{UlmerR2}*{4.7}.  
\item It is evident that if $k$ has characteristic $p$, then starting
  with any $\CC=\DD$ and $f$, $g$ with $g(y)=1/f(y)$ we will arrive at
  an $X$ with unbounded Mordell-Weil rank in the tower of fields $K_d$
  as $d$ varies through the ``special'' values $p^n+1$.  It seems
  that these examples are not all covered by the methods of
  \cite{UlmerR2}.
% and in particular that supersingularity does not play a key role.  
%This and related issues will be pursued further in 
%\cite{OcchipintiThesis}.
\item Still assuming $k$ has positive characteristic, a different
  construction leads to examples of ranks which grow linearly in the
  tower $K_d$.  In other words, we have non-trivial lower bounds on
  ranks for all $d$, not just those dividing $p^n+1$ for some $n$.
  This is one of the main results of \cite{OcchipintiThesis}.
  \end{enumerate}
\end{rems}

\section{Explicit points}\label{s:ExplicitPoints}
The interesting term in the rank formula~\ref{thm:main} is
$\Hom_{k-av}(J_{\CC_d},J_{\CC_d})$ which is a quotient of
$NS(\CC_d\times_k\DD_d)$.  In the example of the previous section,
large ranks came from endomorphisms $\Fr_q\circ[\zeta_d^i]$ of
$J_{\CC_d}$ which in turn are induced by certain explicit maps from
$\CC_d$ to itself.  As we will see in this section, tracing through
the geometry of Berger's construction leads to remarkable explicit
expressions for points on the elliptic curve in the last section.

\begin{thm}\label{thm:points}
  Fix a prime number $p$, let $k=\Fpbar$, and let $K=k(t)$.  Let $X$
  be the elliptic curve 
$$y^2+xy+ty=x^3+tx^2.$$
For each $d$ prime to $p$, let $K_d=k(t^{1/d})\cong k(u)$.  We write
$\zeta_d$ for a fixed primitive $d$-th root of unity in $k$.
\begin{enumerate}
\item For all $d$, the torsion subgroup $X(K_d)_{tor}$ of $X(K_d)$ is
  isomorphic to $\Z/4\Z$ and is generated by $Q=(0,0)$.  We have
  $2Q=(-t,0)$ and $3Q=(0,-t)$.
\item If $p=2$, and $d=p^n+1$, let $q=p^n$ and
$$P(u)=\left(u^q(u^q-u),u^{2q}(\sum_{j=1}^nu^{2^j})\right).$$
Then the points $P_i=P(\zeta_d^iu)$ for $i=0,\dots,d-1$ lie in
$X(K_d)$ and they generate a finite index subgroup of $X(K_d)$, which
has rank $d-1$.  The relation among them is that $\sum_{i=0}^{d-1}P_i$
is torsion.
\item If $p>2$, and $d=p^n+1$, let $q=p^n$ and 
$$P(u)=\left(\frac{u^q(u^q-u)}{(1+4u)^q},
  \frac{u^{2q}(1+2u+2u^q)}{2(1+4u)^{(3q-1)/2}}
  -\frac{u^{2q}}{2(1+4u)^{q-1}}\right).$$ 
Then the points $P_i=P(\zeta_d^iu)$ for $i=0,\dots,d-1$ lie in
$X(K_d)$ and they generate a finite index subgroup of $X(K_d)$, which
has rank $d-2$.  The relations among them are that
$\sum_{i=0}^{d-1}P_i$ and $\sum_{i=0}^{d-1}(-1)^iP_i$ are torsion.
\end{enumerate}
\end{thm}

\begin{rems}\label{rem:explicit}\mbox{}
\begin{enumerate}
\item The expression for $P(u)$ when $p$ is odd is remarkable:  it is
  in closed form (no ellipses or summations) and
  depends on $p$ only through the exponents.
\item The lack of dependence of $P$ on $p$ fits very well with some
  recent speculations of de Jong (cf.~\cite{deJongp09}).
\item We saw in Theorem~\ref{thm:ex1} that $X$ has large rank already
  over $\Fp(t^{1/d})$.  We can find generators for a subgroup of
  finite index by summing the explicit points over orbits under
  $\gal(\Fpbar/\Fp)$.  But it is not clear whether it is possible to
  find expressions for these points which are as explicit and compact
  as those above.
\item We will explain how to get the points in the theorem from
  Berger's construction, but once they are in hand, that they are on
  the curve can be verified by direct calculation and that they
  generate a large rank subgroup can be checked by simple height
  computations which we will sketch in Subsection~\ref{ss:explicit-heights}
  below.
\item As an aid to the reader who wants to check that $P(u)\in
  X(K_d)$, we note that for $p$ odd, the right and left hand sides of
  the Weierstrass equation evaluated at $P(u)$ are both equal to
  $u^{4q}(u^{2q}-2u^{q+1}+u^2)/(1+4u)^{3q-1}$.
\item Ricardo Concei\c c\~ao has recently produced isotrivial elliptic
  curves over $\Fq(t)$ of large rank with many explicit, independent
  polynomial points. See \cite{ConceicaoThesis}.
\end{enumerate}
\end{rems}

\subsection{Plan}
We will sketch the proof of the theorem in the rest of this section.
The calculation of the torsion subgroup is straightforward using the
Shioda-Tate isomorphism $MW(J_X)\cong L^1NS(\XX)/L^2NS(\XX)$ explained
in Subsection~\ref{ss:S-T} and the explicit construction of $\XX_d$
given in Section~\ref{s:Formula}.  We omit the details.

\subsection{Graphs of endomorphisms}
For the rest of the section, we assume $d$ has the form $d=p^n+1$ and
we let $q=p^n$.

Recall that $\en_{k-av}(J_{\CC_d})$ and
$\Hom_{k-av}(J_{\CC_d},J_{\DD_d})$ are quotients of
$NS(\CC_d\times_k\CC_d)$ and $NS(\CC_d\times_k\DD_d)$ respectively,
via the homomorphism which sends a divisor to the action of the
induced correspondence on the Jacobians.  The endomorphism
$\Fr_q\circ[\zeta_d^i]$ of $J_{\CC_d}$ is induced by the graph of the
morphism $\Fr_q\circ[\zeta_d^i]:\CC_d\to\CC_d$.  The statements about
the relations among the $P_i$ follow from Lemma~\ref{lemma:endos}.  To
find explicit coordinates, it suffices to treat the case $i=0$ since
the other points $P_i$ are the images of $P_0$ under $\gal(K_d/K)$.
So from now on we take $i=0$.

Using the isomorphism $\CC_d\to\DD_d$ which sends $(x,z)$ to
$(y,w)=(1/x,1/z)$, the graph of $\Fr_q$ maps to the curve
$$\{y=x^{-q},w=z^{-q}\}\subset
\CC_d\times_k\DD_d$$ 
where $x$,$z$, $y$, and $w$ are the coordinates on $\CC_d$ and $\DD_d$
introduced in Section~\ref{s:eDPC}.

Taking the image in $(\CC_d\times_k\DD_d)/\mu_d$ and using the
isomorphism $(\CC_d^o\times_k\DD_d^o)/\mu_d\cong\XX_d^o$, the graph
maps to the curve
$$\{(x,y,u)=(x,x^{-q},f(x))\}\subset\XX_d^o.$$
(Here $u$ is the coordinate on the $\P^1_k$ in the lower left of the
diagram~\ref{eq:fiber-product} and $f(x)=x(x-1)$.)

Since $f$ has degree 2, the Zariski closure of this curve is a
bi-section of $\pi_d:\XX_d\to\P^1$.  Passing to the generic fiber, we
get a point defined over a quadratic extension.  More precisely, let
$k(a)$ be the quadratic extension of $K_d=k(u)$ with $a(a-1)=u$.  The
point in question has coordinates $(x,y)=(a,a^{-q})$ in the model
$Y=V(f-tg)\subset\CC_d\times_k\DD_d$.  Passing to the Weierstrass
model $X$, our point becomes
$$(x,y)=(-a(1-a)^{q+1},a^2(1-a)^{2q+1}).$$
The next step is to take the trace of this point down to $X(K_d)$

\subsection{Taking a trace}
We use a prime to denote the action of $\gal(k(a)/K_d)$ so that
$a'=(1-a)=b$ and we have $a+b=1$ and $ab=-u$.

We have to compute the sum of
$$(x,y)=(-ab^{q+1},a^2b^{2q+1})$$ 
and 
$$(x',y')=(-a^{q+1}b,a^{2q+1}b^2)$$
in $X(k(a))$, the result being $P=P(u)\in X(K_d)$.  Actually, it will
be more convenient to let $P$ be $-((x,y)+(x',y'))$ (the third point
of intersection of $X$ and the line joining $(x,y)$ and $(x',y')$).

To that end, let
$$\lambda=\frac{y'-y}{x'-x}=\frac{a^{2q}b-ab^{2q}}{(b-a)^q}$$
and
$$\nu=\frac{x'y-xy'}{x'-x}=\frac{(ab)^{q+1}(a^qb-ab^q)}{(b-a)^q}.$$
(Both of these expressions lie in $K_d$ but it is convenient to leave
them in this form for the moment.)

The coordinates of $P$ are then $x_P=\lambda^2+\lambda-u^d-x-x'$ and
$y_P=\lambda x_P+\nu$.

\subsection{The case $p=2$}
Assume that $p=2$.  Then $b-a=1$, and since $a^2=u+a$, $b^2=u+b$,
%\begin{equation*}
\begin{align*}
\lambda&=a^{2q}b-ab^{2q}\\
&=u^q+a^q-a\cr
&=u^q+(a^q-a^{q/2})+\cdots+(a^2-a)\\
&=u^q+u^{q/2}+\cdots+u.
\end{align*}
%\end{equation*}
Similarly $\nu=u^{q+1}(u^{q/2}+\cdots+u)$.  Straightforward
calculation then leads to the coordinates of $P$.

\subsection{The case $p>2$}
Now assume that $p>2$.  Introduce $c=b-a=1-2a$ so that $c'=-c$ and
$c^2=1+4u$.  Then again using $a^2=u+a$ and $b^2=u+b$,
\begin{align*}
\lambda&=\frac{u^q(b-a)+a^q-a}{(b-a)^q}\\
&=\frac{2cu^q+c-c^q}{2c^q}
\end{align*}
and 
\begin{align*}
\nu&=\frac{u^{q+1}(a^q-a)}{(b-a)^q}\\
&=\frac{u^{q+1}(c-c^q)}{2c^q}.
\end{align*}
Straightforward calculation then leads to expressions for $x_P$ and
$y_P$ which are visibly in $k(c^2)=k(u)=K_d$.

This completes our sketch of the proof of the theorem.
\qed

\subsection{Heights}\label{ss:explicit-heights}
Recall that the Mordell-Weil group $X(K_d)$ has a canonical height
pairing $\langle\cdot,\cdot\rangle$ discussed in
Subsection~\ref{ss:heights} above.  Because the points in
Theorem~\ref{thm:points} are so explicit, we can easily compute their
heights.

\begin{prop}\label{prop:heights}
  Let the notation and hypotheses be as in Theorem~\ref{thm:points}.
  We view the points $P_i$ as being indexed by $i\in\Z/d\Z$.
\begin{enumerate}
\item If $p=2$, then
$$\langle P_i,P_j\rangle=
\begin{cases}
\frac{(d-1)^2}{d}&\text{if $i=j$}\\
\\
\frac{(1-d)}{d}&\text{if $i\neq j$}
\end{cases}
$$
\item If $p>2$, then 
$$\langle P_i,P_j\rangle=
\begin{cases}
\frac{(d-1)(d-2)}{d}&\text{if $i=j$}\\
\\
\frac{2(1-d)}{d}&\text{if $i-j$ is even and $\neq0$}\\
\\
0&\text{if $i-j$ is odd}
\end{cases}
$$
\end{enumerate}
\end{prop}

\begin{rems}\mbox{}
\begin{enumerate}
\item The proposition gives another proof that the sums $\sum_iP_i$
  and $\sum_i(-1)^iP_i$ are torsion when $p$ is odd.
\item In the notation of \cite{ConwaySloane}*{4.6.6}, the Mordell-Weil
  lattice $X(K_d)$ modulo torsion is a scaling of $A_{d-1}^*$ when $p=2$
  and a scaling of $A^*_{d/2-1}\oplus A^*_{d/2-1}$ when $p>2$.
\item Let $V_d$ be the subgroup of $E(K_d)$ generated by $Q$ and the
  $P_i$, $i=0,\dots,d-1$ ($V$ for ``visible'').  The BSD formula
  together with the height calculation above yields a simple
  relationship between the order of the Tate-Shafarevitch group of $X$
  over $K_d$ and the index of $V_d$ in $X(K_d)$:
$$|\sha(X/K_d)|=[X(K_d):V_d]^2.$$
\end{enumerate}
\end{rems}

\subsection{Proof of Proposition~\ref{prop:heights}}
First note that $\gal(K_d/K)$ acts on $X(K_d)$ permuting the $P_i$
cyclically and so it will suffice to compute $\langle P_0,P_j\rangle$.
We first treat the case where $p$ is odd.  

We write $P_j$ both for the point in $X(K_d)$ and for the
corresponding section of $\pi_d:\XX_d\to\P^1_k$; similarly for $O$,
the origin of the group law on $X$.  In order to compute $\langle
P_0,P_j\rangle$, we have to find a divisor $D$ supported on components
of fibers such that $P_0-O+D$ is orthogonal to all components of
fibers; then $\langle P_0,P_j\rangle=-(P_0-O+D).(P_j-O)$.  The key
input to finding $D$ and computing the intersection is to find the
reduction of the $P_j$ at each place of $\P^1_k$.

Recall that $P_j$ is the fiberwise sum of a bisection of $\pi_d$ which
we denote $B_j$.  The divisors $B_j-2O$ and $P_j-O$ are linearly
equivalent on $\XX_d$ up to the addition of a sum of components of
fibers, thus in the height calculation we may use $B_j-2O$ in place of
$P_j-O$.  This is a considerable simplification because $B_j$ has a
very compact expression: it is the closure in $\XX_d$ of the point of
$X$ defined over the quadratic extension $k(a)/k(u)$ where
$\zeta_d^ju=a(a-1)$ with coordinates
$$\left(-a(1-a)^{q+1},a^2(1-a)^{2q+1}\right).$$
(With respect to Remark~\ref{rem:explicit} (4), we note that
using $B_j$ is not essential---the heights can be computed knowing
only the $P_j$, not the $B_j$.)

The canonical divisor of $\XX_d$ is $(d/2-2)$ times a fiber.  Since
$B_j$ and $O$ are smooth rational curves, the adjunction formula
implies that $O^2=-d/2$ and $B_j^2=2-d$.  It turns out that $O.B_j=0$
for all $j$, $B_0.B_j=1$ if $j$ is even and $\neq0$, and $B_0.B_j=0$
if $j$ is odd.  So the geometric part of the height pairing is
$$-(B_0-2O).(B_j-2O)=
\begin{cases}
3d-2&\text{if $j=0$}\\
2d-1&\text{if $j$ is even and $\neq0$}\\
2d&\text{if $j$ is odd.}
\end{cases}$$

To find the ``correction factors'' $D.(B_j-2O)$ we calculate the
reduction of $B_j$ at the bad places.  Over $u=0$, $\XX$ has a fiber
of type $I_{4d}$ and if we label the components as usual (with
elements of $\Z/4d\Z$ such that the identity component is 0), the
$B_j$ reduce to distinct points on the components labelled 1 and
$2d+1$ (or $-1$ and $2d-1$ depending on the orientation of the
labelling).  It follows that the correction factor at $u=0$ is
$-(d+2-1/d)$ for all $j$.  Over $u=\infty$, $\XX$ has a fiber of type
$I_d$.  Using the standard labelling of components, if $j$ is even,
$B_j$ reduces to the point of intersection of the components labelled
$d/2-1$ and $d/2$ and if $j$ is odd, $B_j$ reduces to the intersection
of the components labelled $d/2$ and $d/2+1$ (again up to a suitable
choice of orientation).  It follows that the correction factor at
$u=\infty$ is $-(d-1-1/d)$ if $j$ is even and $-(d-2+1/d)$ if $j$ is
odd. Adding the geometric and correction factors gives the result.

The calculations for $p=2$ are quite similar.  One interesting twist
is that at $\infty$, the reduction is of type $I_d^*$ and $B_j$ lands
on one of the components of multiplicity 2.  Thus the linear algebra
required to find the correction term is not carried out in the
standard references.

This completes the proof of the proposition.
\qed

\section{Second example}\label{s:Example2}

\subsection{Data}
Throughout this section, $k=\C$, the field of complex numbers.
Let $\CC=\DD=\P^1_k$.  Let
$$f(x)=x(x-1)(x-a)\quad\text{and}\quad g(y)=y(y-1)(y-a)$$ 
where $a$ is a parameter in $U=k\setminus\{0,1,-1, 1/2,
2,\zeta_6,\zeta_6^{-1}\}$.  (The reason we exclude $-1$, $1/2$, $2$,
and the primitive 6th roots is that $f$ and $g$ have extra symmetry in
those cases and they turn out to be less interesting.) For all $d$, we
have $e_{d,f}=e_{d,g}=1$.

\subsection{The surface $\XX_1$}
The surface $\XX_1$ is obtained by blowing up $\CC\times_k\DD$ once at
each of ten base points.  Straightforward calculation shows that
$\XX_1\to\P^1_k$ is smooth away from $t=0$, $1$, and $\infty$ and two
other points specified below.  The fibers over $0$ and $\infty$
consist of a rational curve of multiplicity 3 and three rational
curves of multiplicity 1 each meeting the multiplicity 3 component at
one point.  (This is a Kodaira fiber of type $IV$ blown up at the
triple point.)  The fiber over $t=1$ is two rational curves of
multiplicity 1 meeting transversely at two points, i.e., a Kodaira
fiber of type $I_2$.  The two other singular fibers are of type $I_1$.
It follows that $c_1(d)=d$ and
$$c_2(d)=
\begin{cases}
4&\text{if $3\nodiv d$}\\
6&\text{if $3| d$.}
\end{cases}$$

\subsection{The curve $X$}
The curve $X$, the smooth proper model of
$\{f-tg\}\subset\P^1_K\times_K\P^1_K$, is an elliptic curve with
invariants
$$c_4=16(a^2-a+1)^2t^2,$$
$$c_6=-216a^2(a-1)^2t^4-16(a-2)^2(a+1)^2(2a-1)^2t^3-216a^2(a-1)^2t^2,$$
and
\begin{equation*}
\Delta=a^2(a-1)^2t^4(t-1)^2Q
\end{equation*}
where
$$Q=-27a^2(a-1)^2t^2
   +(-16a^6 + 48a^5 - 42a^4 + 4a^3 - 42a^2 + 48a - 16)t
   -27a^2(a-1)^2.$$
%  checks for use in Pari:
% c4 = 16*(a^2-a+1)^2*t^2
% c6 =
%   -216*a^2*(a-1)^2*t^4-16*(a-2)^2*(a+1)^2*(2*a-1)^2*t^3-216*a^2*(a-1)^2*t^2
% Delta = a^2*(a-1)^2*t^4*(t-1)^2
%  times -27*a^2*(a-1)^2*t^2+(-16*a^6+48*a^5-42*a^4+4*a^3-42*a^2+48*a-16)*t
%             -27*a^2*(a-1)^2
%
The discriminant of the quadratic $Q$ is 
$$64(a-2)^2(a+1)^2(2a-1)^2(a^2-a+1)^3$$ 
% 64*(a-2)^2*(a+1)^2*(2*a-1)^2*(a^2-a+1)^3
and so for $a\in U$,
$\Delta$ has simple zeroes at the roots of the quadratic.  These roots
are the ``two other points'' referred to above.

\subsection{Homomorphisms and endomorphisms}
The curve $\CC_d$ has genus $d-1$ if $3\nodiv d$ and $d-2$ if $3|d$.
The new part of its Jacobian has dimension $\phi(d)$ if $d>3$ and 1 if
$d=2$ or $3$.

Since $\CC_d\cong\DD_d$ compatibly with the $\mu_d$ action, we have
$$\Hom_{k-av}(J_{\CC_d},J_{\DD_d})^{\mu_d}\cong
\en_{k-av}(J_{\CC_d})^{\mu_d}.$$
Thus our rank formula reads
$$\rk X(K_d)=\rk\en_{k-av}(J_{\CC_d})^{\mu_d}
-d+
\begin{cases}
4&\text{if $3\nodiv d$}\\
6&\text{if $3|d$.}
\end{cases}$$

Since the endomorphisms in $\mu_d$ commute with one another,
Lemma~\ref{lemma:endos} shows that the
minimum value of $\rk\en(J_{\CC_d})^{\mu_d}$ is $d-1$.
The upshot of the following result is that it
is sometimes larger.

\begin{thm}
  For $d\in\{2,5,7,8,9,12,14,15,16,18,22,24\}$, there is a countably
  infinite subset $S_d\subset U$ \textup{(}everywhere dense in the
  classical topology\textup{)} such that if $a\in S_d$,
$$\rk X(K_d)\ge\phi(d)+
\begin{cases}
3&\text{if $3\nodiv d$}\\
5&\text{if $3|d$}
\end{cases}$$
\end{thm}

The largest rank guaranteed by the theorem occurs for $d=15$, $22$,
and $24$ where we have rank $\ge13$.

\begin{proof}
  For any $d>1$ and any $a\in U$ we write $\CC_{d,a}$ for the curve
  with affine equation $z^d=x(x-1)(x-a)$ and $A_{d,a}$ for the new
  part of the Jacobian of $\CC_{d,a}$.  The dimension of $A_{d,a}$ is
  $\phi(d)$ if $d\neq3$ and $1$ if $d=3$.  By Lemma~\ref{lemma:endos},
  the endomorphism algebra of $A_{d,a}$ contains the field
  $\Q(\mu_d)$.  We claim that for $d$ as in the theorem, there are
  infinitely many $a\in U$ such that $A_{d,a}$ is of CM type, more
  precisely such that $\en^0(A_{d,a})$ contains a commutative
  subalgebra of rank $2$ over $\Q(\mu_d)$.  For such $a$,
  $\en(J_{\CC_{d,a}})^{\mu_d}$ has dimension at least $\phi(d)+d-1$.
  Using the rank formula \ref{thm:main} and the calculation of $c_1$
  and $c_2$ above gives the desired result.

  It remains to prove the claim.  The cases $d=5$ and $7$ were treated
  by de Jong and Noot \cite{deJongNoot91}. (See also \cite{deJongAWS}
  for a nice exposition of this and related material.)  The other
  cases are quite similar and so we will just briefly sketch the
  argument.

  First, we may construct $A_{d,a}$ in a family over $U$.  In other
  words, there is an abelian scheme $\sigma:\AF_d\to U$ whose fiber
  over $a\in U$ is $A_{d,a}$.  This abelian scheme gives rise to a
  polarized variation of Hodge structures (PVHS) of weight 1 and rank
  $2\phi(d)$ with endomorphisms by $\Q(\mu_d)$.  (Here and below we
  exclude $d=3$ where the rank is 2 and the variation is locally
  constant.)

  On the other hand, the corresponding period domain for such Hodge
  structures is the product of copies of the upper half plane.  The
  dimension is $r=\{i\,|\,(d,i)=1, (d-1)/3<i<d/2\}$.  For the $d$ in
  the statement, $r=1$.

  At any given point of $U$ we can choose a trivialization of the
  local system $R^1\sigma_*\Z$ in a disk around the point and obtain a
  local period mapping.  A Kodaira-Spencer calculation, which we omit,
  shows that these period maps are not constant.  Since the period
  domain is 1-dimensional, the images of the period maps must contain
  non-empty open subsets.  But the set of CM points is an everywhere
  dense countable subset, and the corresponding values of $a$ are
  values for which the rank of $X_{d,a}$ jumps to at least the value
  stated in the theorem.  This establishes the lower bound stated for
  a countable dense set of values of $a$
\end{proof}

\begin{rems}\mbox{}
\begin{enumerate}
\item The points asserted to exist in this theorem are much less
  explicit than those of Theorem~\ref{thm:points}.  Indeed, even the
  values of $a$ for which the rank jumps are not explicit.
\item The Andre-Oort conjecture would imply that for a fixed large
  $d$, the number of $a$ for which the rank jumps as in the theorem is
  finite.
\item A more interesting question is whether there are arbitrarily large
  $d$ such that there exists at least one $a$ for which the rank
  jumps.  If so, we would have unbounded ranks over $\C(t)$.  It seems
  likely to the author that this does not happen, in other words, that
  for all sufficiently large $d$ and all $a\in U$,
  $\rk\en(A_{d,a})^{\mu_d}=\phi(d)$.  
\item There are at least two other families of curves with countable
  dense CM subsets that could be used in conjunction with the rank
  formula.  See \cite{Desrousseaux04}*{p.~114}
\end{enumerate}
\end{rems}

\end{document}